\documentclass{jhrs-meine}
\usepackage{amsmath}
\usepackage[all]{xy}
\newcommand{\labto}[1]{\xrightarrow{\labelstyle\textstyle #1}}

\usepackage{bbm}
\usepackage[english]{babel}
\usepackage[all]{xy}
\usepackage{amsmath,amssymb,amsthm}
\newtheorem{theorem}{Theorem}
\newtheorem{corollary}{Corollary}
\newtheorem{proposition}{Proposition}
\newtheorem{lemma}{Lemma}
\newtheorem{definition}{Definition}
\def\leq{\leqslant}
\def\geq{\geqslant}
\def\Span{\operatorname{Span}}
\def\Ideal{\operatorname{Ideal}}

\begin{document}                                                        
\title{Homology  Functors  With Cubical Bars}                           
\author{VOLKER W. TH\"UREY}                                             
\email{\\ volker@thuerey.de}                                            
\address{Germany\\28199 Bremen\\Rheinstr. 91\\T: 49 (0)421/591777}      
\begin{abstract}
This work arose from efforts to generalise the usual cubical boundary by using 
different `weights' for opposite faces, but still to obtain a chain complex, 
and this method was found to generalise.
We describe a variant  of the classical singular cubical  homology
theory, in which the usual boundary $(n-1)$-cubes of each $n$-cube
are replaced by combinations of internal $(n-1)$-cubes parallel to
the boundary. This defines a generalised homology theory, but the
usual singular homology can be recovered by taking the quotient  by
the degenerate singular cubes.
\end{abstract}
\startpage{1}                     
\maketitle

\tableofcontents
\newpage
\section{Introduction}
There are different ways to define singular homology groups, for
instance by using simplices, see e.g. \cite{Hatcher}, \cite{Vick},
\cite{Rotman}, or cubes, see \cite{hilton-wylie,Massey}. Because the
latter construction is only one of some possible ways to get this
well known theory, it seems that today mathematicians only have
historical interest in it, but less mathematical interest, because
by using simplices instead of cubes one gets  isomorphic homology
groups, and the simplicial homology theory as it is introduced in
\cite{Eil-Stee} is well-understood and  a common tool of
topologists. Singular homology theory  is a very useful and
successful method not only for  mathematicians but also in other
fields of  science. It is used  for instance  for digital image
processing and nonlinear dynamics, where even the cubical variant is
used, see \cite{Ka-Mi-Mro}. Cubical methods are also essential in
\cite{BHS}.

Here we show an easy way  to generalise cubical singular homology.
In the ordinary  cubical singular homology theory the boundary
operator is constructed by taking the topological boundary of an
$n$-dimensional unit cube as a linear combination of  $2 \cdot n$
cubes of dimension $(n-1)$, provided with alternating signs. We
generalise this by `drawing' in all $n$  directions a linear
combination of a fixed number $L+1$ of   $(n-1)$-dimensional cubes
parallel to the topological boundary, provided  with a coefficient
tuple $ \vec{m} := \left( m_{0}, m_{1}, m_{2},   \dots   ,m_{L}
\right)$. Note that for a fixed  $L > 1$, `our' boundary operator
$_{\vec{m}}\partial_{n}$   is determined  not only by the 
topological boundary but also by parts of the interior of the unit
cube, in contrast to the classical cubical homology theory.

Let  $ \mathsf{TOP}^{2}$ be the  category  of pairs of topological
spaces and continuous maps   as morphisms. That means   that  $ (f:
(X,A) {\rightarrow} (Y,B))  \in \mathsf{TOP}^{2}$ if and only if $X$
and $Y$ are topological spaces  and  $A \subset X,  B \subset Y$ and
$A, B$ carry the subspace topology   and  $f$ is continuous  and
$f(A) \subset B$. Let  ${\cal R}$  be a commutative ring with unit
$1_{\cal R}$. Let  ${\cal R}$-$\mathsf{MOD}$   be the  category of
${\cal R}$-modules.

As in the classical theory, our construction yields a chain complex
with decreasing dimensions, i.e. we get  a sequence of natural
transformations  $ (_{\vec{m}}\partial_{n})_{ n \geq 0 } $ with the
property $_{\vec{m}}\partial_{n}  \circ \, _{\vec{m}}\partial_{n+1}
= 0 $. Hence we shall be able to define homology modules
\[ {}_{\vec{m}}{\cal H}_{n}(X,A) := \frac {\text{kernel}(_{\vec{m}}\partial_{n })}
{\text{image}(_{\vec{m}}\partial_{n+1} )},\] and this  will lead  to
a  sequence of functors   $_{\vec{m}}{\cal H}_{n}:  \mathsf{TOP}^{2}
\longrightarrow {\cal R}$-$\mathsf{MOD}$,  for all $n \geq 0$.

The exactness axiom follows immediately, and with an additional
condition on the fixed coefficient tuple $ \vec{m}$ the homotopy
axiom holds. Unfortunately, so far the excision axiom could be
verified only in the case of $L=1$,  but in that special case we get
a class of extraordinary homology theories.  For $L=1$, our boundary
operator can be regarded as a kind of `weighted' topological
boundary of a cube, with the weight   $ \left( m_{0}, m_{1} \right)
$.

In this way for every  fixed   $L \in \mathbbm{N}$   and  fixed
tuple $\vec{m} \in  {\cal R} ^{L+1}$  a functor $_{\vec{m}}{\cal
H}_{n} : \mathsf{TOP}^{2} \longrightarrow  {\cal R}$-$\mathsf{MOD}$
will  be  constructed  for each  $n  \geq 0$. In the special case
${\cal R} :=   \mathbbm{Z} $   we shall see that the   homotopy
axiom  holds if and only if the   greatest common divisor of
$\{m_{0}, m_{1}, m_{2}, \dots  ,m_{L} \}$   is  $1$. If  $L = 1$ and
gcd \{$ m_{0}, m_{1}$\} = 1    the excision axiom holds, so we get
an extraordinary  homology theory.    Finally, the ordinary singular
homology can be recovered by taking quotients by `degenerate'
singular cubes.

This new construction is primarily of theoretical interest, since
all the homology modules which we can compute in our homology theory
we can already compute in terms of the ordinary  singular homology,
by using the clever method from \cite{BuCoFl}.    If we take the
ring   ${\cal R} :=  \mathbbm{Z} $, we are able to compute the
homology groups   for a  finite {\it CW}-complex, and we express
them as   a product of   singular homology groups.

We assume that the reader is familiar both with the construction of
the classical cubical singular homology, e.g. in \cite{Massey}, and
also with homological algebra and ordinary singular homology theory,
see   e.g.  \cite[p.57 ff]{Rotman}, or \cite[p.97 ff]{Hatcher}.

\section{General Definitions and Notations}
We denote  the natural numbers by  $\mathbbm{N}_0 := \{0, 1, 2, 3,
\dots \}$, the positive integers by $ \mathbbm{N} :=  \{1, 2, 3,
\dots \} $, the ring of integers by $ \mathbbm{Z}$ and the real
numbers by   $\mathbbm{R}$.

The brackets    $\left( \cdots \right)$  will be used for tuples and
besides  $\left[ \cdots \right]$ to  structure  text and formulas,
$\left[ r , s \right]$    also for the closed interval, [{\sf u}]
for the equivalence  class of a quotient module. The brackets
$\left\langle \cdots  \right\rangle$    will be needed for the
boundary operator, \ $\left\| \cdots \right\|$  for the subdivision
operator and  $\left\{ \cdots \right\}$  for sets.

As before let ${\cal R}$ be a commutative ring with unit $1_{\cal
R}$.   Let $X$  be a topological space. All maps we shall use will
be continuous.

For each $ n \in \mathbbm{N}$ let  ${ \bf I }^{n}$ be the {\it
$n$-dimensional unit cube}, that means ${ \bf I }^{n}  :=  \left\{
\,  \left( x_{1}, \; x_{2}, \; \ldots  , \;x_{n} \right) \in
\mathbbm{R}^{n} \ | \  x_{i}   \in   \left[ 0 ,1 \right] \ \text{for
} \ 1 \leq i \leq n \, \right\} $, provided with the Euclidean
topology, and let ${ \bf I }^{0}  :=  \{ 0 \}$. We write ${ \bf I
}^{1} = { \bf I } = [ 0, 1 ] $,   the unit interval.
\begin{definition}  \rm  \quad   \label{definition eins}    Define the sets
${\cal S}_n(X)   :=   \left\{ T : { \bf I }^{n} \rightarrow  X \ | \
T \ \text{is continuous} \; \right\}$,     and $ {\cal K}_{  n}(X)
:= $  the free ${\cal R}$-module with the basis  ${\cal S}_n(X)$,
for all   $n \in \mathbbm{N}_0 $, as well as  $  {\cal K}_{  -1}(X)
:= \{ 0 \} $, the trivial  ${\cal R}$-module. Every   {\sf u} $\in
{\cal K}_{  n}(X) $   is called a {\it chain}.  $ { } $   \hfill
$\Box$
 \end{definition}
Let us assume given for each topological space $X$ and $n \in
\mathbbm{N}_0$ an ${\cal R}$-module morphism $\partial_{n} \colon
{\cal K}_{n}(X) \longrightarrow  {\cal K}_{n-1}(X)$. If we have the
property   $ \partial_{n} \circ \partial_{n+1} = 0$ for all $n \geq
0$,  we call the map  $\partial_{n}$  a {\it boundary operator}, and
the sequence $(\partial_{n } )_{n \geq 0}$ of ${\cal R}$-module
morphisms is called a {\it chain complex }  $ {\cal K}_*(X) $,
\[ {\cal K}_*(X) \ := \ \cdots \cdots \
\labto{\partial_{n+1}} {\cal K}_{n}(X)     \labto{\partial_{n}}
{\cal K}_{n-1}(X) \labto{\partial_{n-1}}   \ \cdots \
\labto{\partial_1} {\cal K}_{0}(X) \labto{\partial_0} \{ 0 \}   .
   \]
An element {\sf u} $\in \text{ kernel} (_{\vec{m}}\partial_{n})$  is
called a  {\it cycle}, an element  {\sf w} $\in \text{ image}
(_{\vec{m}}\partial_{n+1})$  is called a  {\it boundary}. Because
$\partial_{n} \circ \, \partial_{n+1 } =  0 $   the   ${\cal
R}$-module
\begin{align*}
{\cal H}_{n}(X) := \frac {kernel (\partial_{n })}{image
(\partial_{n+1})}
\end{align*}
is well defined for all  topological spaces $X$ and all   $ n \in
\mathbbm{N}_0$;  $ {\cal H}_{n}(X)$   is called the  {\it
$n^{th}$-homology ${\cal R}$-module of $X$}.

Two continuous functions  $f , g \colon X  \longrightarrow  Y$  are
{\it homotopic}, written $ f  \simeq  g$,  if  and  only  if there
is a  continuous $ H : X \times { \bf I }  \longrightarrow Y $ such
that for all   $x \in X$ we have   $ H (x,0) = f(x)$  and  $H(x,1) =
g(x)$.
\begin{definition}    \label{definition zwei}  \rm   \quad
Let ${\cal H} := ({\cal H}_{n})_{n \geq 0 }  $   be a sequence of
functors $ {\cal H}_{n}:    \mathsf{TOP} \longrightarrow    {\cal
R}$-$\mathsf{MOD} $. We say that  the functors $ {\cal H}$   satisfy
the {\it Homotopy Axiom} if  and  only  if for all  $ f , g :  X
\longrightarrow  Y$  with  $f  \simeq  g $   we have ${\cal
H}_{n}(f)  =  {\cal H}_{n}(g)$, for all  $n \in  \mathbbm{N}_0 $.
\hfill  $\Box$
\end{definition}
Let  ${\cal M}$ be an arbitrary subset of the ring ${\cal R}$. We
write $\Span _{\cal R}({\cal M})$ for the ideal of  ${\cal R}$
generated by  ${\cal M}$. We have  $\Span _{\cal R}({\cal M})  =
{\cal R}$   if and only if there are a number $k \in \mathbbm{N}$
and sets $\{ m_{1},  m_{2}, \ldots , m_{k}\} \subset {\cal M}$ and $
\{ r_{1},  r_{2}, \ldots ,  r_{k} \} \subset  {\cal R}$ such that
$\sum_{i = 1}^{k}  r_{i} \cdot m_{i}  =  1_{\cal R}$.        \hfill
$\Box$
\begin{definition}    \rm  \quad  \label{definition drei}
Let  $ a, b $ be elements of the ring  ${\cal R}$. We say that
$(a,b)$  fulfils the condition ${\cal {NCD}}$   if and only if for
all    $n   \in   \mathbbm{N}$  there exists $ x_n ,   y_n \in
{\cal R}$  with
\begin{equation*}
x_n \cdot a^{n} +  y_n   \cdot b^{n} = 1_{\cal R} .
\tag*{$\Box$}\end{equation*}
\end{definition}
Of course, ${\cal {NCD}}$ is equivalent to $\Span _{\cal R}({a^{n},
b^{n}})  = {\cal R}$  for all  $n  \in  \mathbbm{N}$.   The letters
${\cal {NCD}}$   remind us of   `\! No Common Divisor'.   In the
ring  $ \mathbbm{Z}$ we have that  $ (a,b) $ has the property ${\cal
{NCD}}$  if and only if the ideal generated  by  $\{a , b\}$ is
$\mathbbm{Z}$.
\begin{definition}    \rm   \quad    \label{definition vier}
Let $\Im$  be a set of indices, let ${\cal U}$ := $\{ U_i  \mid i
\in \Im \}$
be a family of subsets of  $X$ whose interiors cover  $X$. Let  \\
\centerline{   $  {\cal S}_{n}(X,\: {\cal U}) := \left\{ T \in {\cal
S}_n(X)  \: | \: \text{there is an}  \:
i  \in  \Im \ \text{such that} \ T({ \bf I }^{n})  \subset U_i \right\}$.    }     \\
For all    $n  \in   \mathbbm{N}_{0}$    and for every topological
space $X$ we define ${\cal K}_{ n}(X,{\cal U})$ to be  the  free
${\cal R}$-module  with  the  basis   $  {\cal S}_{n}(X,\: {\cal U})
$. The elements   {\sf u} $\in  {\cal K}_{ n}(X,{\cal U}) \; $ are
called {\it ${\cal U}$-small  chains}.  For a subset $A  \subset  X$
with the canonical inclusion $ i: A  \hookrightarrow  X$, the map
$i$ leads to  a canonical inclusion
$$\widehat{i}: {\cal K}_{ n}(A,{\cal U}) \; \hookrightarrow {\cal K}_{ n}(X,{\cal U})$$
in   ${\cal R}$-$\mathsf{MOD}$.  Define  $${\cal K}_{ n}( X,A,{\cal
U} )  := \frac{{\cal K}_{ n}( X,\: {\cal U} )}{{\cal K}_{ n}( A,\:
{\cal U} )},$$ and  this yields an inclusion
$${\cal K}_{ n}( X,A,{\cal U} ) \stackrel{j}{\hookrightarrow} {\cal K}_{ n}( X,A )$$ \hfill  $\Box$
\end{definition}
\begin{definition}   \rm  \quad    \label{definition fuenf}
For each $ T \in {\cal S}_n(X)$, i.e.  $T : { \bf I }^{n}
\rightarrow  X $,  we construct the maps $\left\| T \right\|_{\alpha
, \vec{e} , \vec{v} } \in {\cal S}_n(X)$, which will be used for the
excision axiom. First we need an auxiliary map $q_n :
\mathbbm{R}^{n} \rightarrow  { \bf I }^{n}$, for all dimensions $ n
\in \mathbbm{N}$. For  $( y_1 , y_2 , \ldots , y_{n} ) \in
\mathbbm{R}^{n} $ let $q_n( y_1 , y_2 , \ldots , y_{n} )   :=   (
z_1 , z_2 , \ldots , z_{n} )\in  {\bf I }^{n}$, where
$$   z_i   :=  \begin{cases} 0   & \text{  if  }       y_i  \leq    0   ,    \\
y_i & \text{  if  }       y_i  \in    [0,1]  ,  \\
1   & \text{  if  }       y_i  \geq    1  , \quad
 \text{ for }  i \in \{ 1 , 2, \ldots ,  n \}.
 \end{cases} $$
For fixed  $ \alpha  \in  \mathbb{R}$  and fixed  $\vec{v} := ( v_1
, v_2 , \ldots , v_{n} ),    \vec{e} := ( e_1 , e_2 , \ldots , e_{n}
) \in  \mathbb{R}^{n}$, define the map $H_{\alpha , \vec{e} ,
\vec{v} } : { \bf I }^{n} \rightarrow  \mathbb{R}^{n}$ by $H_{\alpha
, \vec{e} , \vec{v} } \: ( x_1 , x_2 , \ldots , x_{n} ) := ( y_1 ,
y_2 , \ldots , y_{n} )$,     where
$$ y_i := \alpha \cdot ( e_i + v_i \cdot x_i ), \quad \text{for all }  i \in \{ 1, 2, \ldots , n \}.  $$
Finally for each    $T : { \bf I }^{n} \rightarrow  X $ we set
$\left\| T \right\|_{\alpha , \vec{e} , \vec{v} }   := (T \circ q_n
\circ H_{\alpha , \vec{e} , \vec{v} } )\colon {\mathbf I}^n \to  X
$. \hfill  $\Box$
\end{definition}
\begin{definition} \rm   \label{definition sechs}
Let  $ n \in  \mathbbm{N}$.  Let $ {\cal  D}_{n}(X) $ be the  subset
of $ {\cal S}_{n}(X)$ consisting of the  {\it degenerate } cubes,
i.e. those $ T \in {\cal S}_{n}(X)$ such that  there is a  $j \in
\{1,2, \ldots , n\}$ and for $ y , z   \in   [ 0,1 ]$  \ we have
$$ T ( x_1,  \ldots , x_{j-1} , y , x_{j+1} , \ldots , x_n )  =
T ( x_1, \ldots , x_{j-1} , z , x_{j+1} , \ldots , x_n ) ,$$ i.e.
$T$ does not depend on the $j^{\: th}$  component.
 We shall write  \\
         \centerline{   $  T ( x_1, x_2, \ldots , x_{j-1} , * , x_{j+1} , \ldots , x_n ) \ := \
                  T ( x_1, x_2, \ldots , x_{j-1} , y , x_{j+1} , \ldots , x_n )$. }  $ { } $  \hfill  $\Box$
     \end{definition}
     \begin{definition}     \rm  \quad    \label{definition sieben}
           We define the free  ${\cal R}$-module ${\cal K}_{{\cal D},  n}(X)$, which will be an ideal of
                  ${\cal K}_{ n}(X)$.
                  A chain   $\sum_{i=1}^{p} r_i \cdot T_i \in  {\cal K}_{ n}(X)$
                  is an element of ${\cal K}_{{\cal D},  n}(X)$  if and only if for all
                  $i = 1, 2, \ldots, p$   we have  $T_i   \in  {\cal  D}_{n}(X) $.
                  That means  that ${\cal K}_{{\cal D},  n}(X)$ is the free ${\cal R}$-module generated
                  by degenerate maps.  We have ${\cal K}_{{\cal D},0}(X) = \{ 0 \} $.  \hfill  $\Box$
     \end{definition}
     \begin{definition}    \rm \quad    \label{definition acht}
            Define   for a fixed    ${\alpha} \in {\cal R}$ the  submodule $\Ideal _{\alpha,  n}(X)$  of
                  ${\cal K}_{ n}(X)$,  generated by ${\alpha} {\cal R}$.
                  That means  a chain $\sum_{i=1}^{p} r_i \cdot T_i  \in  {\cal K}_{ n}(X)$  belongs to
                 $\Ideal _{\alpha,  n}(X)$   if and only if   for all   $ i = 1, 2, \ldots, p$
                  \, there is an element $y_i \in {\cal R}$   such that    $r_i = y_i \cdot \alpha$.
                   $ { } $   \hfill  $\Box$
     \end{definition}
     \begin{definition}     \rm     \label{definition neun}
              For fixed $\alpha \in {\cal R}$ and  for all  $n \in
              \mathbbm{N}_0$ let
         $$ \Gamma_{\alpha,  n }(X) :=\Ideal _{\alpha,  n}(X)   +
                            {\cal K}_{{\cal D},  n}(X)$$      (generally this is not a direct sum).
             Then $ \Gamma_{\alpha,   n }(X)$   is an ideal of   ${\cal K}_{ n}(X)$,   and a chain
              $\sum_{i=1}^{p} r_i \cdot T_i  \in  {\cal K}_{ n}(X)$  belongs to $ \Gamma_{\alpha,   n }(X)$
              if and only if   for all   $ i = 1, 2, \ldots, p$  either $ T_i $   is degenerate, or
              $ r_i$ is a multiple of $ \alpha $.        \hfill $\Box$
    \end{definition}
    \begin{definition}      \rm  \quad    \label{definition zehn}
              Correspondingly to the previous three definitions  we define  for pairs  $ (X,A) \in  \mathsf{TOP}^{2}$
              for a fixed $ \alpha \in {\cal R} $   the quotients of
              ${\cal R}$-modules $ {\cal K}_{{\cal D},  n}(X,A)$, $\Ideal _{\alpha,  n}(X,A)$   and
              $ \Gamma_{\alpha, n }(X,A)$,  for all  $ n \in  \mathbbm{N}_0$.
              That means that two chains  {\sf u} and  { \sf w}  represent the same
              equivalence class if and only if  the difference $ {\sf u} - {\sf w}$ is a chain in   $A$.
              Further, we define  the quotients of ${\cal R}$-modules
             $${\cal K}_{ n}(X)_{\sim  \Gamma, \alpha }  :=
                 \frac{{\cal K}_{ n}(X)}{\Gamma_{\alpha,   n }(X)} \quad \ \text{and} \quad \
                 {\cal K}_{ n}(X,A)_{\sim  \Gamma, \alpha }  :=
                 \frac{{\cal K}_{ n}(X,A)}{\Gamma_{\alpha,   n }(X,A)} \ \  .    $$   \hfill  $\Box$
    \end{definition}

  \section{The Boundary Operator}      \label{section boundary operator}
           Fix a  natural  number $L \geq 1 $  and an $(L+1)$-tuple  $\vec{m}$ of ring elements,
             $\vec{m} := \left( m_{0}, m_{1}, m_{2},   \dots   ,m_{L} \right)  \in  {\cal R}^{L+1}.$
             For each  $ n \in  \mathbbm{N}_0$  we shall define  a `boundary operator'
             $_{\vec{m}}\partial_{n}:  {\cal K}_{  n}(X) \stackrel{ }{\longrightarrow} {\cal K}_{n-1}(X)$.
      For integers  $ n \geq 1$ and  $T \in {\cal S}_n(X) $ we shall define a map
             $\left\langle T\right\rangle_{ n, i, j } \in  {\cal S}_{n-1}(X)$,
             for all integers $ 0 \leq i \leq L $ and $ 1 \leq j \leq  n $.

        For  $n = 1$  let  $\left\langle T\right\rangle_{ 1, i, 1 }   \in   {\cal S}_{0}(X) $
             with  $\left\langle T\right\rangle_{ 1, i, 1 }  (0)  := T\left( \frac{i}{L}\right)$. \
             For $ n > 1$ and every ($n-1$)-tuple $( x_{1}, x_{2}, \dots , x_{n - 1} )\ \in {\bf I }^{n - 1}$ we set
       $$\left\langle T\right\rangle_{ n, i, j }  ( x_{1}, x_{2}, \dots , x_{j - 1}, x_{j}, \dots , x_{n - 1} )\
              :=  T \left( x_{1}, x_{2}, \dots , x_{j - 1}, \frac{i}{L}, x_{j}, \dots , x_{n - 1 } \right). $$
              Finally for every $n  \in  \mathbbm{N}$ and all $T \in {\cal S}_{n}(X) $ let
      \begin{align}
                _{\vec{m}}\partial_{n}(T) :=  \sum_{j =  1}^{n} (-1)^{j+1} \cdot   \sum_{i = 0}^{L}
                                              m_{i} \cdot \left\langle T\right\rangle_{ n, i, j } ,
      \end{align}
      and for $ n = 0$ let $_{\vec{m}}\partial_{0}(T) := 0$, the only possible map.

      See Figure \ref{picture1}, which illustrates the case   $n := 2,   L := 3,  {\vec{m}} := (9,1,4,-3)$   and
               $ T := id({ \bf I }^{2})$. On the left hand side you see the two-dimensional unit cube
               ${ \bf I }^{2}$, the right hand side shows   $ _{\vec{m}}\partial_2 (T)$, i.e. the images of
               eight one-dimensional unit cubes                                           
               $  \left\langle T\right\rangle_{ 2, i, j } , i \in \{0,1,2,3\}$ and $j \in
               \{1,2\} $, multiplied by  coefficients   $ 9,1,4,-3 $, elements of the ring ${\cal R} := \mathbbm{Z}$.
  \begin{figure}[ht]
           \setlength{\unitlength}{1cm}
           \begin{picture}(8,5)                      
                 \put(0.4,0){\line(1,0){3.0}}
                 \put(0.4,1){\line(1,0){3.0}}
                 \put(0.4,2){\line(1,0){3.0}}
                 \put(0.4,3){\line(1,0){3.0}}
                 \put(0.4,-0.0){\line(0,1){3.0}}      \put(0.1,-0.4){ $ 0 $}
                 \put(1.4,-0.0){\line(0,1){3.0}}      \put(1.17,-0.4){ $ \frac{1}{3} $}
                 \put(2.4,-0.0){\line(0,1){3.0}}      \put(2.18,-0.4){ $ \frac{2}{3} $}
                 \put(3.4,-0.0){\line(0,1){3.0}}      \put(3.18,-0.4){ $ 1 $}
                 \put(2.9,0){\vector(1,0){1}}         \put(3.775,-0.4){x$_1$}
                 \put(0.4,2.5){\vector(0,1){0.8}}     \put(0.24,3.5){x$_2$}
                 \put(0,0.9){ $ \frac{1}{3} $}
                 \put(0,1.9){ $ \frac{2}{3} $}
                 \put(0,2.8){ $ 1 $}
                 \put(9.6,0){\line(1,0){3.0}}         \put(5.05,1.4){9}
                 \put(9.6,1){\line(1,0){3.0}}         \put(5.6,1.4){+ 1}
                 \put(9.6,2){\line(1,0){3.0}}         \put(6.6,1.4){+ 4}
                 \put(9.6,3){\line(1,0){3.0}}         \put(7.45,1.4){ $-$ 3}
                 \put(5.4,0){\line(0,1){3.0}}         \put(8.8,-0.1){$-$ 9}
                 \put(6.4,0){\line(0,1){3.0}}         \put(8.8, 0.9){$-$ 1}
                 \put(7.4,0){\line(0,1){3.0}}         \put(8.8, 1.9){$-$ 4}
                 \put(8.4,0){\line(0,1){3.0}}         \put(8.8,2.9){+ 3}
         \end{picture}
 \caption{} \label{picture1}
 \end{figure}
           \\    \\     \\  
                 For a `chain'   {\sf u} = $r_{1} \cdot T_{1} + r_{2} \cdot T_{2} \in  {\cal K}_{n}(X)$
                   define   $_{\vec{m}}\partial_{n}$({\sf u})  by linearity, \\
                 \centerline{  $_{\vec{m}}\partial_{n}(r_{1} \cdot T_{1} + r_{2} \cdot T_{2}) := r_{1} \cdot \,                                    _{\vec{m}}\partial_{n} ( T_{1} ) + r_{2} \cdot \, _{\vec{m}}\partial_{n} ( T_{2}$). } \\
        { \bf Remark}:  For    $L   \geq   2$    the map    $_{\vec{m}}\partial_{n}$
                  is  defined  not only on  the topological boundary  but also on parts of the interior of
                  ${ \bf I }^{n}$.  We use the name `boundary operator' only  for historical reasons.
      \begin{theorem}    \label{erstes Theorem}
            For arbitrary  $L \in \mathbbm{N}$ and $(L+1)$-tuples
            $\vec{m} = \left( m_{0}, m_{1}, m_{2}, \dots  ,m_{L} \right) \in  {\cal R}^{L+1} $
            we have for all $n  \in \mathbbm{N}: \ _{\vec{m}}\partial_{n - 1} \; \circ \; _{\vec{m}}\partial_{n } = 0 $.
      \end{theorem}
    \begin{proof} \quad
           For   $n = 1$   the statement is trivial.
                   For  $n = 2$  the proof is similar to  the cases with  $n \geq 3$, which we shall show in  detail.
                   Thus, let  $n \geq 3$.
                   Because of the linearity of   $_{\vec{m}}\partial_{n }$   it suffices  to prove the theorem
                   for the basis  of  $ {\cal K}_{  n}(X) $.  Let  $T : { \bf I }^{n} \stackrel{}{\rightarrow} X$
                   be continuous, i.e.  $ T \in {\cal S}_n(X)$.  We have
            \begin{align}
                  _{\vec{m}}\partial_{n - 1} \; \circ \;    _{\vec{m}}\partial_{n }  (T)
                  &  =  \,  _{\vec{m}}\partial_{\:n - 1}
                         \left(    \sum_{j =  1}^{n}   (-1)^{j+1} \cdot   \sum_{i = 0}^{L}
                         m_{i} \cdot  \left\langle T\right\rangle_{ n, i, j }   \right)     \\
                  &  =   \sum_{j =  1}^{n}   (-1)^{j+1}  \cdot   \sum_{i = 0}^{L}   m_{i} \cdot \,
                         _{\vec{m}}\partial_{n - 1} \left( \left\langle T\right\rangle_{n, i, j} \right) \\
                  &  =   \sum_{j =  1}^{n}   (-1)^{j+1}  \cdot   \sum_{i = 0}^{L}
                              m_{i}  \cdot  \sum_{p =  1}^{n-1}   (-1)^{p+1}  \sum_{k = 0}^{L}  m_{k} \cdot
                              \left\langle  \left\langle T\right\rangle_{n, i, j}
                              \right\rangle_{n-1, k, p}.\notag \\
    \intertext{Thus}                      \label{equation drei}
                       _{\vec{m}}\partial_{n - 1} \; \circ \; _{\vec{m}}\partial_{n } (T)
                &=  \sum_{j =  1}^{n}   \sum_{p =  1}^{n-1}   (-1)^{j+p+2} \cdot \sum_{i = 0}^{L}
               \sum_{k = 0}^{L}  m_{i}  \cdot   m_{k}  \cdot
               \left\langle  \left\langle T\right\rangle_{n, i, j} \right\rangle_{n-1, k, p}.
         \end{align}
            Note that the following display is not suited for special cases like $ j=1, j=n, p=1, p=n-1,j=p$
            or $ j=p+1$.  Nevertheless the claim of Theorem  \ref{erstes Theorem} remains true in all these cases.
            Let  $( x_{1} ,\  x_{2} ,\  x_{3} ,\  \dots \; ,  x_{n - 2} ) $ be an arbitrary point in
            ${ \bf I }^{n - 2}$. For fixed  $i, k  \in \{ 0, 1, 2, \dots ,  L \}$   we have for
            $ p  \in \{ 1, 2, \dots , n-1 \}$ and  $ j  \in \{ 1, 2, \dots , n \}$
       \begin{align*}
              &\left\langle  \left\langle T\right\rangle_{ n, i, j }
                     \right\rangle_{n-1, k, p} (x_{1}, \ldots , x_{p-1}, x_{p},  \ldots , x_{n - 2})   \\
              &   =   \left\langle T\right\rangle_{n,  i, j } \left( x_{1}, \ldots , x_{p-1} , \frac{k}{L} ,
                       x_{p} ,  \ldots ,  x_{n-2}  \right)   \\
              &   =
                  \begin{cases}
                       T \left(  x_{1} , \ldots , x_{j-1} , \frac{i}{L} , x_{j}, \ldots , x_{p-1} , \frac{k}{L} ,
                        x_{p} , \ldots  , x_{n-2}  \right)   & \text{ if } j \leq  p  ,     \\
                       T  \left(  x_{1} ,   \ldots , x_{p-1} , \frac{k}{L} , x_{p}, \ldots , x_{j-2} , \frac{i}{L}  ,
                        x_{j-1} , \ldots , x_{n-2} \right) & \text{ if }j >  p  .
                  \end{cases}
         \end{align*}
   The sign of $_{\vec{m}}\partial_{n - 1}  \circ \, _{\vec{m}}\partial_{n }(T) $   depends on
                $j$ and $p$  only. It is easy to see that for  $j \leq  p$  we get
  \begin{equation*}
    (-1)^{j+p+2} \,   m_{i} \, m_{k} \left\langle  \left\langle T\right\rangle_{ n, i, j  }
                \right\rangle_{\: n-1, k, p } +  \ (-1)^{(p+1)+j+2} \,  m_{k} \, m_{i} \left\langle
                \left\langle T\right\rangle_{ n, k, p+1 } \right\rangle_{ n-1,  i, j }  = 0 .
  \end{equation*}
               The set    $M :=  \left\{ 1, 2, 3, \dots  ,n \right\}  \times
              \left\{ 1, 2, 3, \dots  ,n-1 \right\}$  contains  $n \cdot (n-1)$  elements.
              With     $M_{small}  :=  \{ (j,p) \in M | j \leq p \}$
              and      $M_{big}  :=  \{ (j,p) \in M | j > p \}$    we have
              $M = M_{small}  \cup  M_{big}$,  and  $M_{small} \cap M_{big} = \emptyset$.
              The map \\
      \centerline{ $M_{small} {\rightarrow} M_{big}, \ (j,p)   \mapsto   (p+1,j)$  }  \\
              is  bijective. Thus  the  $ n \cdot (n-1) \cdot (L+1)^{2}$  maps in equation (\ref{equation drei}) 
              cancel pairwise. Hence
              $_{\vec{m}}\partial_{n - 1} \;     \circ \; _{\vec{m}}\partial_{n }(T) = 0 $.
 \end{proof}
     For all  topological spaces  $X$ the homology groups
     $ _{\vec{m}}{\cal H}_{n}(X) = \frac {kernel ( _{\vec{m}}\partial_{n })}{image ( _{\vec{m}}\partial_{n+1 })}$
     of the chain complex \ $_{\vec{m}}{\cal K}_*(X) \ = $
   \[                      \cdots \cdots   \  \labto{ _{\vec{m}}\partial_{n+1}} {\cal K}_{n}(X)
                                                \labto{ _{\vec{m}}\partial_{n}}   {\cal K}_{n-1}(X)
                                                \labto{ _{\vec{m}}\partial_{n-1}}    \ \cdots \
               \labto{ _{\vec{m}}\partial_1} {\cal K}_{0}(X)  \labto{ _{\vec{m}}\partial_0} \{ 0 \}
   \]
         are well defined since $_{\vec{m}}\partial_{n} \circ \, _{\vec{m}}\partial_{n+1} = 0 $,
         for all fixed  $\vec{m}  \in {\cal R}^{L+1},   n \in \mathbbm{N}_0$.
         For a cycle {\sf u}, i.e.  $_{\vec{m}}\partial_{n}({\sf u}) = 0$,  we denote the  equivalence  class
         containing {\sf u}  by [{\sf u}]$_{\sim}   \in \, _{\vec{m}}{\cal H}_{n}(X)$.

            We use the abbreviation  $ _{\vec{m}}{\cal H} := (_{\vec{m}}{\cal H}_n)_{n \geq 0}$.
            We say that the fixed tuple $ \vec{m} \in {\cal R}^{L+1} $ is the
            {\it weight},  the element  $\sigma := \sum_{i =  0}^{L} \; m_{i} \in {\cal R}$
            is the  {\it index} of $ _{\vec{m}}{\cal H}$. The number $L \in \mathbbm{N}$
            is called the {\it  length} of the weight $\vec{m}$.

     Example: \quad  For the one-point space $\{p\}$ and for  $n \in \mathbbm{N}_0$  there is only one
                    $T: { \bf I }^{n} \rightarrow \{p\}$, thus we have ${\cal K}_{n}({p}) \cong  {\cal R}$.                                And for the chain complex $_{\vec{m}}{\cal K}_*(p)$,
    \[ _{\vec{m}}{\cal K}_*(p) \ = \ \cdots  \cdots \
                   \labto{_{\vec{m}}\partial_4} {\cal K}_{3}(p)    \labto{_{\vec{m}}\partial_3} {\cal K}_{2}(p)
                   \labto{_{\vec{m}}\partial_2} {\cal K}_{1}(p)    \labto{_{\vec{m}}\partial_1} {\cal K}_{0}(p)
                   \labto{_{\vec{m}}\partial_0} \{ 0 \}   ,
    \]
           we  get  \quad  $ _{\vec{m}}{\cal K}_*(p) \ \cong \ \cdots  \cdots \
                   \labto{_{\vec{m}}\partial_4}   {\cal R}
                   \labto{_{\vec{m}}\partial_3}   {\cal R}            \labto{_{\vec{m}}\partial_2} {\cal R}
                   \labto{_{\vec{m}}\partial_1}   {\cal R}            \labto{_{\vec{m}}\partial_0} \{ 0 \}$.    \\
           If we define the map  $\times \sigma: {\cal R} \rightarrow {\cal R}$,
                $ x \mapsto \sigma \cdot x$, we can describe the  boundary operators  by
        \[    _{\vec{m}}\partial_{n}
                   \cong \begin{cases} 0               & \text{  if  }   n   \text{ is even } , \\
                                       \times \sigma   & \text{  if  }   n   \text{ is odd}.
                         \end{cases}
        \]
            Explanation: Note the definition of   $_{\vec{m}}\partial_{n}(T) = \sum_{j =  1}^{n}   (-1)^{j+1} \cdot
                     \sum_{i = 0}^{L} m_{i} \cdot \left\langle T\right\rangle_{ n, i, j }$ .   Because of
                     the alternating signs  $\sigma$  copies of the unique map from   ${\bf I}^{n-1}$ to $\{p\}$
                     cancel pairwise. That means for an arbitrary index  $\sigma $   that
      \[  _{\vec{m}}{\cal H}_{n}(p)
                     \cong \begin{cases} \{ x \in  {\cal R} \, | \, \sigma \cdot x = 0 \}
                                                                & \text{  if  } n  \text{ is odd }  \\
                     {\cal R}/(\sigma \cdot {\cal R})           & \text{  if  } n  \text{ is even}  ,
                     \end{cases}
      \]
        i. e. for  $\sigma = 0$   we get
              $ _{\vec{m}}{\cal H}_{n}(p)   \cong   {\cal R}  $  for all $n \in  \mathbbm{N}_0$.     \\

       The above  construction of  $_{\vec{m}}{\cal H}_{n}(X) $ yields a functor
              $ _{\vec{m}}{\cal H}_{n}:  \mathsf{TOP} \longrightarrow {\cal R}$-$\mathsf{MOD}$, for all
              $n \in  \mathbbm{N}_0$:  Let  $f: X \rightarrow Y$  be  continuous  and
              $T \in   {\cal S}_{n}(X) $, then we have $ f \circ T  \in  {\cal S}_{n}(Y) $.
              Let ${\cal K}_{  n}(f):  {\cal S}_{  n}(X)  \longrightarrow  {\cal S}_{  n}(Y) $,
              for all  basis elements   $ T \in  {\cal S}_{n}(X) $   let
              $ {\cal K}_{  n}(f)(T)  :=  f \circ T $.  We create a functor
              ${\cal K}_{  n} : \mathsf{TOP} \longrightarrow {\cal R}$-$\mathsf{MOD}$,  with \\
       \centerline { $  {\cal K}_{  n} \left( X \stackrel{ f }{\longrightarrow} Y \right) :=
                 {\cal K}_{  n}(X)  \xrightarrow{ {\cal K}_{  n}(f) }   {\cal K}_{  n}(Y) $,   }   \\
        $  {\cal K}_{  n}(f)$  is  well defined by linearity,
               and for an arbitrary    $( f : X  \rightarrow Y ) \in   \mathsf{TOP}$    the following
               diagram commutes  in      ${\cal R}$-$\mathsf{MOD}$     for all  $n  \in \mathbbm{N}$:
    \[   \xymatrix{   \cdots \quad    \ar[r]^<>(0.4){_{\vec{m}}\partial_{n+2}}
            & {\cal K}_{ n+1}(X)  \ar[d]<0.0cm>^<>(0.45){{\cal K}_{  n+1}(f)}
                                        \ar[r]^<>(0.5){_{\vec{m}}\partial_{n+1}}
          &  {\cal K}_{ n}(X)   \ar[d]<0.0cm>^<>(0.45){{\cal K}_{  n}(f)}
                                      \ar[r]^<>(0.5){_{\vec{m}}\partial_{n}}
          &  {\cal K}_{ n-1}(X) \ar[d]<0.0cm>^<>(0.45){{\cal K}_{  n-1}(f)}
                                      \ar[r]^<(0.4)>(0.5){_{\vec{m}}\partial_{n-1}}
          &  {\cal K}_{ n-2}(X) \ar[d]<0.0cm>^<>(0.45){{\cal K}_{  n-2}(f)}
                                      \ar[r]^<>(0.5){_{\vec{m}}\partial_{n-2}}
          &                        \quad        \cdots           \\
          \cdots   \quad  \ar[r]_<>(0.4){_{\vec{m}}\partial_{n+2}}
            &  {\cal K}_{ n+1}(Y)  \ar[r]_<>(0.5){_{\vec{m}}\partial_{n+1}}
          &  {\cal K}_{ n}(Y)  \ar[r]_<>(0.5){_{\vec{m}}\partial_{n}}
          &  {\cal K}_{ n-1}(Y) \ar[r]_<(0.0)>(0.5){_{\vec{m}}\partial_{n-1}}
          &  {\cal K}_{ n-2}(Y)   \ar[r]_<>(0.5){_{\vec{m}}\partial_{n-2}}
          &        \quad   \cdots      \\     }
    \]
 Thus we get   ${\cal K}_{  n-1}(f) \circ \, _{\vec{m}}\partial_{n }  =  \,
 _{\vec{m}} \partial_{n } \circ {\cal K}_{  n}(f) $   for all   $n \in  \mathbbm{N}_0$,
 hence    $ {\cal K}_{  n}(f)$   maps  cycles to cycles and  boundaries to boundaries.
 For an arbitrary map $( f : X  \longrightarrow Y ) \in \mathsf{TOP}$, for a cycle {\sf u},
 hence  [{\sf u}]$_{\sim}   \in \, _{\vec{m}}{\cal H}_{n}(X)$,    let \\
 \centerline{  $_{\vec{m}}{\cal H}_{n}(f)$ ([{\sf u}]$_{\sim}  )   :=
                                 [{\cal K}_{n}(f)(${\sf u}$)]_{\sim} \in \, _{\vec{m}}{\cal H}_{n}(Y)$,    } \\
 and  we define  \\
 \centerline{  $_{\vec{m}}{\cal H}_{n} \left(  X  \stackrel{f}{\longrightarrow} Y \right)   := \
 _{\vec{m}}{\cal H}_{n}(X)  \xrightarrow{ _{\vec{m}}{\cal H}_{n}(f) } \
 _{\vec{m}}{\cal H}_{n}(Y)$.    }  \\
 In this way   $_{\vec{m}}{\cal H}_{n} $   is a functor
 $  \mathsf{TOP} \longrightarrow   {\cal R}$-$\mathsf{MOD}$.

    In a similar way $_{\vec{m}}{\cal H}_{n}$   will be extended to a functor
                      $  \mathsf{TOP}^{2} \longrightarrow {\cal R}$-$\mathsf{MOD}$:
             (The following description is rather brief. For more details the reader should  study
             \cite[p.95 ff]{Rotman},  or \cite[p.22 ff]{Massey}, or other  books  about singular homology theory).

    If there is  $( f: (X,A) \rightarrow (Y,B) )  \in   \mathsf{TOP}^{2}$,
             we have  subspaces  $A   \hookrightarrow X$ and $B \hookrightarrow Y$  (in $ \mathsf{TOP}$),
             and submodules  $ {\cal K}_{  n}(A) \hookrightarrow  {\cal K}_{  n}(X) $  as well as
             $ {\cal K}_{  n}(B) \hookrightarrow {\cal K}_{  n}(Y)$  (in   ${\cal R}$-$\mathsf{MOD}$).
             Hence the   ${\cal R}$-modules
             $  {\cal K}_{  n}(X,A) :=  \frac{ {\cal K}_{  n}(X) }{ {\cal K}_{  n}(A) }$
             and     $ {\cal K}_{  n}(Y,B) :=  \frac{ {\cal K}_{  n}(Y) }{ {\cal K}_{  n}(B) }$
             are well defined for   $n  \in \mathbbm{N}_0 \,  , \; {\cal K}_{  -1}(X,A) := \{0\} $.
             For  a  `chain'   {\sf u} $\in  {\cal K}_{  n}(X) $  let  [{\sf u}] $\in  {\cal K}_{  n}(X,A)$
             be the equivalence  class of  {\sf u}  modulo   $  {\cal K}_{  n}(A)$.
             Thus  [{\sf u}] = [{\sf w}]  if and only if  {\sf u} $-$ {\sf w} $\in   {\cal K}_{  n}(A)$,
             that means that  {\sf u} $-$ {\sf w}  is a chain in  $A$.

The just constructed boundary operator  $ _{\vec{m}}\partial_{n}:
{\cal K}_{  n}(X)  \longrightarrow {\cal K}_{  n-1}(X)  $ also
yields a map $ {\cal K}_{  n}(X,A)$  $\longrightarrow {\cal K}_{
n-1}(X,A)$, which we call  $ _{\vec{m}}\partial_{n} $,  too.  It has
the property $ _{\vec{m}}\partial_{n} \circ \,
_{\vec{m}}\partial_{n+1} = 0 $, for $ n \geq 0$, as before. We
define a corresponding chain complex \, $ _{\vec{m}}{\cal K}_*(X,A)$
for pairs $ (X,A) $, to be
\[  \cdots \cdots \   \labto{_{\vec{m}}\partial_{n+1}}  {\cal K}_{n}(X,A) \labto{ _{\vec{m}}\partial_{n}}
{\cal K}_{n-1}(X,A)\labto{ _{\vec{m}}\partial_{n-1}} \ \cdots \
 \labto{ _{\vec{m}}\partial_1} {\cal K}_{0}(X,A) \labto{ _{\vec{m}}\partial_0} \{ 0 \}  .\]
 For  ($ f: (X,A) \rightarrow (Y,B))  \in   \mathsf{TOP}^{2}$    and    $T \in  {\cal S}_{n}(A) $
 we have   $f \circ T \in  {\cal S}_{n}(B) $  \enlargethispage{0.5cm}
 (because    $ f (A) \subset B$).   Let
 $$ {\cal K}_{n}(f) ([{\sf u}])  := [{\cal K}_{  n}(f) ({\sf u})] \in {\cal K}_{n}(Y,B)
                                           \ \text{ for }  [{\sf u}] \in  {\cal K}_{  n}(X,A) , $$
       and   $ {\cal K}_{n}(f)$  is well defined. Hence ${\cal K}_{n}$   yields  a functor
             $ \mathsf{TOP}  \longrightarrow {\cal R}$-$\mathsf{MOD}$    as well as a functor
             $ \mathsf{TOP}^{2}  \longrightarrow {\cal R}$-$\mathsf{MOD}$.

As above,  for an arbitrary $( f: (X,A) \longrightarrow (Y,B) )  \in
\mathsf{TOP}^{2}$, the following diagram commutes in    ${\cal R}$-$
\mathsf{MOD}$   for    $n  \in \mathbbm{N}_0$:
    \[   \xymatrix{   \cdots \quad    \ar[r]^<>(0.4){_{\vec{m}}\partial_{n+2}}
            &  {\cal K}_{ n+1}(X,A)     \ar[d]<0.0cm>^<>(0.45){{\cal K}_{  n+1}(f)}
                                        \ar[r]^<>(0.5){_{\vec{m}}\partial_{n+1}}
          &  {\cal K}_{ n}(X,A)       \ar[d]<0.0cm>^<>(0.45){{\cal K}_{  n}(f)}
                                      \ar[r]^<>(0.5){_{\vec{m}}\partial_{n}}
          &  {\cal K}_{ n-1}(X,A)     \ar[d]<0.0cm>^<>(0.45){{\cal K}_{  n-1}(f)}
                                      \ar[r]^<(0.4)>(0.5){_{\vec{m}}\partial_{n-1}}
          &                        \quad        \cdots           \\
          \cdots   \quad              \ar[r]_<>(0.4){_{\vec{m}}\partial_{n+2}}
            &  {\cal K}_{ n+1}(Y,B)     \ar[r]_<>(0.5){_{\vec{m}}\partial_{n+1}}
          &  {\cal K}_{ n}(Y,B)       \ar[r]_<>(0.5){_{\vec{m}}\partial_{n}}
          &  {\cal K}_{ n-1}(Y,B)     \ar[r]_<(0.0)>(0.5){_{\vec{m}}\partial_{n-1}}
          &                                                                    \quad   \cdots      \\     }
    \]
    For a chain   {\sf u}  $\in  {\cal K}_{n}(X)$,  by abuse of notation we sometimes denote the
          equivalence class  [{\sf u}] $ \in  {\cal K}_{n}(X,A)$  simply by  {\sf u}.
          If  $_{\vec{m}}\partial_{n }$({\sf u})$ \in {\cal K}_{n-1}(A) $
          let from now on  [{\sf u}]$_{\sim} \in \, _{\vec{m}}{\cal H}_{n}(X,A)$ be the equivalence class  modulo                 $\text{ image} (_{\vec{m}}\partial_{n+1})$.  We call such a chain  {\sf u} a {\it cycle}. For a cycle
          {\sf u} and   $ f : (X,A) \rightarrow (Y,B) $, let   \\
     \centerline{$_{\vec{m}}{\cal H}_{n}(f)$ ([{\sf u} ]$_{\sim}) := [{\cal K}_{n}(f)$({\sf u})]$_{\sim} \, , $ } \\
      and we have  \\
     \centerline{ $ _{\vec{m}}\partial_{n } \circ {\cal K}_{n}(f) ({\sf u}) =  {\cal K}_{n-1}(f)
                    \circ \, _{\vec{m}}\partial_{n} ({\sf u}) \in {\cal K}_{  n-1}(B)$.  }   \\

      Hence   ${\cal K}_{n}(f)$   maps  cycles to cycles and  boundaries to boundaries as above. Therefore
          $\left[ {\cal K}_{n}(f)({\sf u}) \right]_{\sim} \in \, _{\vec{m}}{\cal H}_{n}(Y,B)$,
          and  the functor $_{\vec{m}}{\cal H}_{n}:  \mathsf{TOP}^{2}  \longrightarrow {\cal R}$-$ \mathsf{MOD}$
          is well defined for all  $n \in \mathbbm{N}_0$.    As an abbreviation take   \\
    \centerline { $_{\vec{m}}{\cal H}_{n} \left[ (X,A) \stackrel{f}{\longrightarrow}  (Y,B) \right]   =: \
              _{\vec{m}}{\cal H}_{n}(X,A) \stackrel{f_\ast}{\longrightarrow} \, _{\vec{m}}{\cal H}_{n}(Y,B)$. }

    \section{The Homotopy Axiom}
          The reader should note that in the following we shall omit the weight `${\vec{m}}$' in the boundary operator
          $_{\vec{m}}\partial_n$ for an easier display.
    \begin{lemma}    \label{lemma eins}  \quad
           $_{\vec{m}}{\cal H}$  satisfies the homotopy axiom   if and only if
           for all   topological spaces   $X$    and for   $e_0 , e_1 : X \rightarrow  X \times  { \bf I }$,
           $e_0 (x) := (x,0)$    and    $e_1(x) : = (x,1)$,   the equation
           $ _{\vec{m}}{\cal H}_n(e_0) = \, _{\vec{m}}{\cal H}_n(e_1)$   holds for every $n \in \mathbbm{N}_0$.              \end{lemma}
     \begin{proof} \quad
           `$\Longrightarrow$':  \ \
                Since  $e_0  = Id_{X \times { \bf I }}  \circ  e_0$  and
                $e_1  = Id_{X \times { \bf I }}  \circ  e_1$,  we have
                $e_0  \simeq  e_1$.     \\
           `$\Longleftarrow$':  \ \
                If we assume    $f \simeq g$  there is a continuous $H$ with   $ f = H \circ  e_0$    and \
                $g = H \circ  e_1$,    and   $ _{\vec{m}}{\cal H}_n$   is  a functor, hence    \\                                       \centerline{ $ _{\vec{m}}{\cal H}_n(f) = \ _{\vec{m}}{\cal H}_n(H) \circ \ _{\vec{m}}{\cal H}_n(e_0) =                     \ _{\vec{m}}{\cal H}_n(H) \circ \ _{\vec{m}}{\cal H}_n(e_1) = \  _{\vec{m}}{\cal H}_n(g)$. }
      \end{proof}
      \setlength{\parindent}{0mm}
      Note that the maps   $e_0, e_1 : X \rightarrow  X \times { \bf I }$    induce  canonically  two maps  \\
       \centerline{  $e_0, e_1 :  {\cal S}_n(X)   \rightarrow   {\cal S}_{n}(X \times { \bf I })$, \,
                     by $e_i(T) := e_i \circ T$, for $i \in \{0,1\} $,  }   \\
       and by linearity  two maps $e_0, e_1 : {\cal K}_{ n}(X) \rightarrow {\cal K}_{ n}(X \times {\bf I})$.
  \begin{theorem}[Homotopy Axiom]  \label{Homotopie Axiom}
        Let $L \in \mathbbm{N}$ and  let  $\vec{m} = ( m_{0}, m_{1}, \ldots  ,m_{L}) \in {\cal R}^{L+1}$ be the
        weight of $_{\vec{m}}{\cal H}$. Then $_{\vec{m}}{\cal H}$  satisfies the homotopy axiom if and only if
        $\Span _{\cal R}(\vec{m})  = {\cal R}$.
  \end {theorem}
     \setlength{\parindent}{10mm}
     \begin{proof}  \quad   `$\Longleftarrow$': \quad  We assume that $\Span _{\cal R}(\vec{m})  = {\cal R}$.
                  Because of this assumption a set    $\{ r_{0},  r_{1}, \ldots , r_{L} \}  \subset {\cal R}$
                  exists with  $\sum_{k = 0}^{L}  r_{k} \cdot m_{k}  =  1_{\cal R}$.

         First  we construct  $L+1$  continuous auxiliary functions.  For fixed   $L  \in  \mathbbm{N}$  and
                  for all  $k  \in  \{ 0 , 1 , 2 , \ldots  ,  L \}$  we define  a map
                  $\chi_{k} : [ 0 , 1 ] \rightarrow  [ 0 , 1 ] $.
                  The functions   $\chi_{k}$   are mostly  $ 0 $ and they have  a  `jag'  of  height $1$ at
                   $\frac{k}{L}$ .  More precisely, for  $ k=0$ and $k=L $   the formulas are  \\
        \centerline {  $   \chi_0(x)  :=
                   \begin{cases}  1 - L \cdot x  &  \text{ for } \quad x  \in \left[ 0, \frac{1}{L} \right] , \\
                   0   &  \text{ for } \quad x  \in  \left[ \frac{1}{L}, 1
                   \right],
                   \end{cases} $  }   \\   \\  \\
        and    \centerline {    $
          \chi_L(x)  :=      \begin{cases}  0  &   \text{ for } \quad x  \in  \left[ 0, \frac{L-1}{L} \right] , \\
                             L \cdot x - L + 1 &  \text{ for } \quad x  \in  \left[ \frac{L-1}{L}, 1 \right]  \ .
                             \end{cases}   $ }  \\
               For $ L > 1 $  and $k  \in  \{ 1 , 2, \ldots ,  L-1  \}$  we define
               $\chi_k$  to be  the polygon in $\mathbb{R}^2$ through the five points   (0,0), $\left( \frac{k-1}{L},0 \right)$,
               $\left( \frac{k}{L},1 \right)$, $\left( \frac{k+1}{L},0 \right)$  and
               $(1,0)$, as given by the following formula and
               Figure \ref{pic2} for the case $(k,L)=(2,5)$.

  $$ \qquad  \qquad \qquad  \chi_k(x)        :=
                 \begin{cases} 0  & \text{ \rm for } \quad x  \in
                 \left[ 0, \frac{k-1}{L} \right] \cup  \left[ \frac{k+1}{L}, 1 \right] , \\
                 L \cdot x - k + 1  & \text{ \rm for } \quad x  \in  \left[ \frac{k-1}{L}, \frac{k}{L} \right] , \\
                 1 - L \cdot x + k  & \text{ \rm for } \quad x  \in  \left[ \frac{k}{L}, \frac{k+1}{L} \right]  \ .
                 \end{cases}    $$
  \medskip

  \begin{center}

         \begin{figure}[ht]
                \setlength{\unitlength}{0.5cm}
                \begin{picture}(8,3)    \thinlines
                     \put(0,0){\vector(1,0){7.0}}   \put(6.90,-0.5){$x$}
                     \put(1,-0.5){\vector(0,1){5.5}}   \put(0.5,5.05){$y$}
                      \thicklines     \put(1,0.04){\line(1,0){1.0}}  \put(2,0.03){\line(1,4){1.0}}
                      \put(3,4.03){\line(1,-4){1.0}}    \put(4,0.04){\line(1,0){2.0}}
                      \put(0.55,3.7){$1$}    \put(0.9,4.0){\line(1,0){0.2}}
                      \put(0.55,-0.6){$0$}   \put(1.8,-0.85){$\frac{1}{5}$}   \put(2.8,-0.85){$\frac{2}{5}$}
                      \put(3.8,-0.85){$\frac{3}{5}$}   \put(4.8,-0.85){$\frac{4}{5}$}   \put(5.8,-0.7){$1$}
                      \put(4.9,1.6){The polygon $ \chi_2$,     $ L = 5$}
                      \put(2.0,-0.1){\line(0,1){0.2}}   \put(3.0,-0.1){\line(0,1){0.2}}
                      \put(4.0,-0.1){\line(0,1){0.2}}   \put(5.0,-0.1){\line(0,1){0.2}}
                      \put(6.0,-0.1){\line(0,1){0.2}}
               \end{picture}
  \caption{} \label{pic2}
       \end{figure}
 \end{center}
        Note that for    $j,k  \in  \{ 0 , 1 , 2 , \ldots  ,  L \}$    we have
               $\chi_{k}(\frac{j}{L})   = \delta_{j,k}$ (that means  $\chi_{k}(\frac{j}{L})   = 1 $   if
               $k = j$   and    $\chi_{k}(\frac{j}{L})  = 0 $    if   $k \neq j$).

        Now we  define `chain homotopies',
               $\Theta_n:  {\cal K}_{ n}(X) \;\longrightarrow   {\cal K}_{ n+1}(X \times { \bf I })$,
               which means that the  $\Theta_n$'s will satisfy the equation
               $\partial_{n+1} \circ  \Theta_n = \pm ( e_0  -  e_1 ) +
               \Theta_{n-1}  \circ \partial_{n}$.
               More precisely,  for  $n  \in \mathbbm{N}$  and  $k \in \{ 0 , 1 , 2 , \ldots , L \}$   define
               $\xi_{n}, \psi_{n,k}:  {\cal S}_n(X) \longrightarrow {\cal S}_{n+1}(X \times { \bf I })$
               as  follows.       For every    $T : { \bf I }^{n} \longrightarrow  X$,  for all
               $(x_1 , x_2 , \ldots , x_n, x_{n+1}) \in  { \bf I }^{n+1}$   let
     \begin{align*}
            \xi_{n}(T)   (x_1 , x_2 , \ldots , x_n, x_{n+1}) & := \left( T(x_1, x_2, \ldots , x_{n}), 0 \right),  \\
            \psi_{n,k}(T) (x_1 , x_2 , \ldots ,  x_n, x_{n+1})
            &  :=  \left( T(x_1, x_2, \ldots , x_{n}),  \chi_{k}(x_{n+1}) \right) ,
     \end{align*}
      and  for $n = 0$ let   $\xi_{0}(T) (x) := ( T(0),  0 ) $
                and  $\psi_{0,k}(T) (x) := ( T(0),  \chi_{k}(x) )$, for   $ x \in { \bf I }$.
                Finally let  $\Theta_{-1} := 0$,  and for all  $n  \in \mathbbm{N}_0$ we set
      \begin{align}
                \Theta_n(T) :=  \sum_{k = 0}^{L} r_k  \cdot ( \xi_{n}(T) - \psi_{n,k}(T) ).
      \end{align}
         For    {\sf u}  $\in {\cal K}_{ n}(X)$ let $\Theta_n$({\sf u}) be  defined by
              linearity.  Hence we get for all integers  $ n  \geq -1 $ an ${\cal R}$-linear  map
              $\Theta_n : {\cal K}_{ n}(X) \longrightarrow {\cal K}_{ n+1}(X \times { \bf I })$.
                Thus  we get the following (noncommutative)  diagram  in   ${\cal R}$-$ \mathsf{MOD}$.
    \[   \xymatrix{   \cdots \quad
            \ar[r]^<>(0.5){\partial_{n+2}}
            &  {\cal K}_{ n+1}(X)        \ar[dd]<0.1cm>^{e_1}   \ar[dd]<-0.15cm>_{e_0}
                                       \ar@<0.0 cm>[ldd]_<>(0.4){ \Theta_{n+1}}
                                       \ar[r]^<>(0.5){\partial_{n+1}}
          &  {\cal K}_{ n}(X)          \ar[dd]<0.1cm>^{e_1}  \ar[dd]<-0.15cm>_{e_0}
                                       \ar[ldd]_<>(0.4){ \Theta_{n}} \ar[r]^<>(0.5){\partial_{n}}
          &  {\cal K}_{ n-1}(X)        \ar[dd]<0.1cm>^{e_1}  \ar[dd]<-0.15cm>_{e_0}
                                       \ar[ldd]_<>(0.4){ \Theta_{n-1}}  \ar[r]^<(0.4)>(0.5){\partial_{n-1}}
          &  \quad  \cdots             \ar[ldd]_<>(0.4){ \Theta_{n-2}}               \\
          \cdots   \quad    & & & &    \quad  \cdots                                 \\
          \cdots  \quad
            \ar[r]_<>(0.5){\partial_{n+2}}
            &  {\cal K}_{ n+1}(X \times { \bf I })     \ar[r]_<>(0.5){\partial_{n+1}}
          &  {\cal K}_{ n}(X \times { \bf I })       \ar[r]_<>(0.5){\partial_{n}}
          &  {\cal K}_{ n-1}(X \times { \bf I })     \ar[r]_<(0.0)>(0.5){\partial_{n-1}}
          &      \quad                               \cdots        }
    \]
   \begin{lemma}  \quad  \label{lemma zwei}
        For all  $n \in   \mathbbm{N}_0$ and all $T \in  {\cal S}_n(X) $  we get
        $$  \left[  \partial_{n+1} \circ  \Theta_n \right] (T)  =  \left[ (-1)^{n+2} \cdot
            ( e_0 -  e_1 ) +  \Theta_{n-1}  \circ \,  \partial_{n} \right](T) .  $$
  \end{lemma}
  \begin{proof} \quad
            For  $n = 0$  the proof is a  simpler version of the following one and will be omitted.
            For  $n \in  \mathbbm{N}$  we have:
     \begin{align*}
            \left[ \partial_{n+1} \circ \Theta_n \right] (T)
            &  =    \partial_{n+1}    \left[   \sum_{k=0}^{L}   r_k \cdot
            (  \xi_{n}(T)  -  \psi_{n,k} (T)   )\ \right]    \\
            &  =   \sum_{k=0}^{L}   r_k \cdot \sum_{j=1}^{n+1}   (-1)^{j+1} \cdot \sum_{i=0}^{L}   m_i  \cdot
            \left[   \left\langle  \xi_{n}(T) \right\rangle_{n+1,i,j} - \left\langle
            \psi_{n,k}(T) \right\rangle_{n+1,i,j}     \right]     \\
            &  =  Rest +  {\cal D} , \quad \text{where}
     \end{align*}
     \begin{align}     \label{equation vier}
            Rest & := \ \sum_{j=1}^{n} (-1)^{j+1} \cdot \sum_{i,k=0}^{L}   r_k
            \cdot  m_i  \cdot  \left[ \left\langle  \xi_{n}(T) \right\rangle_{n+1,i,j} - \left\langle  \psi_{n,k}(T)
            \right\rangle_{n+1,i,j}  \right]  , \\   \text{and} \quad
            {\cal D}  & := \ (-1)^{n+2} \cdot \sum_{i,k=0}^{L}   r_k \cdot m_i \cdot
            \left[ \left\langle  \xi_{n}(T) \right\rangle_{n+1,i,n+1} -
                     \left\langle  \psi_{n,k}(T) \right\rangle_{n+1,i,n+1} \right] .
      \end{align}
                 We have  for $ i \in \{ 0,1, \ldots , L \} $  and all tuples
                 $( x_1 , x_2 , \ldots , x_{n} )  \in { \bf I }^{n}$:
       \begin{align*}
                &    \left[ \left\langle  \xi_{n}(T) \right\rangle_{n+1,i,n+1} - \left\langle  \psi_{n,k}(T)
                         \right\rangle_{n+1,i,n+1} \right] ( x_1 , x_2 , \ldots , x_{n} )  \\
                & =  \xi_{n}(T) \left( x_1 , x_2 , \ldots , x_{n},\frac{i}{L} \right)
                           -   \psi_{n,k}(T) \left( x_1 , x_2 , \ldots , x_{n},\frac{i}{L} \right)  \\
                & = ( T( x_1 , x_2 , \ldots , x_{n} ), 0 )   -   \left( T ( x_1 , x_2 , \ldots , x_{n}),
                         \chi_{k}(\frac{i}{L}) \right)  .
       \end{align*}
      Since $\chi_{k}\left(\frac{i}{L}\right) = \delta_{i,k}$ and $\sum_{k=0}^{L} r_k \cdot m_k = 1_{\cal R}$
      it follows that
      \begin{align*}
              {\cal D}  & =  (-1)^{n+2} \cdot \sum_{k=0}^{L}    r_k \cdot  m_k  \cdot
              \left[  \left\langle  \xi_{n}(T) \right\rangle_{n+1,k,n+1} - \left\langle  \psi_{n,k}(T)                                               \right\rangle_{n+1,k,n+1}    \right]  \\                                                                       &  =  (-1)^{n+2} \cdot \sum_{k=0}^{L}    r_k \cdot  m_k   \cdot
                                  \left[      e_0 \circ T  -   e_1 \circ T     \right]  \\
              &  =  (-1)^{n+2} \cdot ( e_0 \circ T  -   e_1 \circ T ) \cdot \sum_{k=0}^{L}    r_k \cdot  m_k  \\
              &  =  (-1)^{n+2} \cdot ( e_0 \circ T - e_1 \circ T ) \ = \ (-1)^{n+2} \cdot ( e_0 (T) - e_1 (T) )  .
      \end{align*}
      It remains  to show that
        $\left[ \Theta_{n-1}\ \circ \, \partial_{n} \right](T) = Rest$.  We   have
     \begin{align*}
                  & \left[ \Theta_{n-1}\ \circ  \partial_{n} \right]  (T)
                     =   \Theta_{n-1}\     \left(  \sum_{j =  1}^{n}   (-1)^{j+1} \cdot \sum_{i = 0}^{L}
                                  m_{i} \cdot  \left\langle T\right\rangle_{ n, i, j } \right)    \\
                  &  =   \sum_{j =  1}^{n}   (-1)^{j+1} \cdot \sum_{i = 0}^{L}  m_{i} \cdot
                                 \Theta_{n-1}\ \left(   \left\langle T\right\rangle_{ n, i, j }  \right) \\                               &  =   \sum_{j =  1}^{n}   (-1)^{j+1} \cdot \sum_{i = 0}^{L}    m_{i} \cdot
                                 \sum_{k = 0}^{L}    r_k  \cdot
                                 \left( \xi_{n-1}(  \left\langle T\right\rangle_{ n, i, j }  ) -
                                 \psi_{n-1,k}(  \left\langle T\right\rangle_{ n, i, j } ) \right) \\
                  &  =   \sum_{j =  1}^{n}   (-1)^{j+1} \cdot \sum_{i , k = 0}^{L}  m_{i} \cdot  r_k \cdot
                                \left(   \xi_{n-1}( \left\langle T\right\rangle_{ n, i, j } ) -
                                 \psi_{n-1,k}( \left\langle T\right\rangle_{ n, i, j }  )  \right)  .                          \end{align*}
       We consider the maps    $ \xi_{n-1}$   and   $\psi_{n-1,k}$   more carefully.
                 We have for all integers   $n > 1   , $     for   $T \in  {\cal S}_n(X) $,    for   $j \in
                 \{1, 2, \ldots , n\}$,   and    $i , k \in \{0, 1, 2, \ldots , L\}$   the following two
                  equations for each   $n$-tuple   $( x_1 , x_2 , \ldots , x_{n-1}, x_{n} ) \in { \bf I }^{n}$:
       \begin{align*}
                 \xi_{n-1} \left( \left\langle T\right\rangle_{ n, i, j } \right)( x_1 , x_2 , \ldots , x_{n})
                 &  =  \left( \left\langle T\right\rangle_{n, i, j }
                                                 ( x_1 , x_2 , \ldots , x_{n-1}),  0 \right) \\
                 &  =  \left(  T ( x_1 , x_2 , \ldots , x_{j-1} , \frac{i}{L} , x_j , \ldots ,
                                                 x_{n-1}),  0 \right)  \\
                 &  =  \xi_{n}(T) \left( x_1 , x_2 , \ldots , x_{j-1} , \frac{i}{L} , x_j , \ldots , x_{n} \right) \\
                 &  =  \left\langle  \xi_{n}(T) \right\rangle_{n+1,i,j}  ( x_1 , x_2 , \ldots , x_{n} )  .
          \end{align*}
          Shortly, we have  $\xi_{n-1} \left( \left\langle T\right\rangle_{ n, i, j } \right)$
             =   $ \left\langle \, \xi_{n}(T) \, \right\rangle_{n+1,i,j}$   .  \\
           In the same way we  find   that
        \begin{align*}
            \psi_{n-1,k} \left( \left\langle T\right\rangle_{ n, i, j } \right) ( x_1 , x_2 , \ldots , x_{n})
            &  =    \left( \left\langle T\right\rangle_{ n, i, j } ( x_1 , x_2 , \ldots , x_{n-1}),                                                                            \chi_{k}(x_{n}) \right) \\
            &  =    \left(  T ( x_1 , x_2 , \ldots , x_{j-1} , \frac{i}{L} , x_j , \ldots ,
                                                                 x_{n-1}), \chi_{k}(x_{n}) \right) \\
            &  =    \psi_{n,k}(T)  \left( x_1 , x_2 , \ldots , x_{j-1} , \frac{i}{L} , x_j , \ldots , x_{n} \right) \\
            &  =    \left\langle \, \psi_{n,k}(T) \, \right\rangle_{n+1,i,j}  ( x_1 , x_2 , \ldots , x_{n} ) ,
       \end{align*}
             therefore  $\psi_{n-1,k} \left( \left\langle T\right\rangle_{ n, i, j } \right)$
             =   $ \left\langle \psi_{n,k}(T) \right\rangle_{n+1,i,j}$,
             and   finally  $\left[ \Theta_{n-1}\ \circ \, \partial_{n} \right] (T) = Rest$, as defined in
                 equation (\ref{equation vier})  follows.
             This ends the proof of Lemma \ref{lemma zwei}.
  \end{proof}
           We have just proved that   $[ \partial_{n+1}\circ \, \Theta_n](T) $
                    =  $[(-1)^{n+2} \cdot ( e_0 -  e_1 ) + \Theta_{n-1} \circ \, \partial_{n}](T) . $
                   Take a cycle {\sf u} $\in  {\cal K}_{ n}(X)$
                   (i.e. $ \partial_{n}$({\sf u}) = 0) instead of $T$.  The fact that
              $$ [ \partial_{n+1} \circ \Theta_n ] ({\sf u}) = [ (-1)^{n+2} ( e_0 -  e_1)] ({\sf u})
                  \  \in \text{ image} (\partial_{n+1})$$
              means that $( e_0 -  e_1)({\sf u})$ is a boundary, hence we can deduce for the equivalence class
                    [{\sf u}]$_{\sim} \in \, _{\vec{m}}{\cal H}_{n}(X)$  that
           $ _{\vec{m}}{\cal H}_{n}( e_0 -  e_1)([{\sf u}]_{\sim})  = 0 =
            \, _{\vec{m}}{\cal H}_{n}(e_0)([{\sf u}]_{\sim}) - \,  _{\vec{m}}{\cal H}_{n}(e_1)([{\sf u}]_{\sim}), $
                    and therefore  $_{\vec{m}}{\cal H}_{n}(e_0) = \, _{\vec{m}}{\cal H}_{n}(e_1)$.
                    By  Lemma \ref {lemma eins}  the  homotopy axiom is satisfied.  \\ \\
    `$\Longrightarrow$': \quad
        We assume that   $_{\vec{m}}{\cal H}$   satisfies the homotopy axiom.   We fix an  $n \in \mathbbm{N}$.
            Let $X := { \bf I }^{2}$,   and let   $T_1 , T_2 :  { \bf I }^{n}  \rightarrow { \bf I }^{2}$
            with $T_1 \neq T_2$, but $ \partial_{n}( T_1 ) = \, \partial_{n}( T_2 ) $. \;  Let  \\
       \centerline{  {\sf u} := $T_1 - T_2 \in {\cal K}_{ n}(X)$. }     \\
        Because $ \partial_{n}$ ({\sf u}) = 0,  {\sf u} is a cycle,  hence
             [{\sf u}]$_\sim \in  \, _{\vec{m}}{\cal H}_{n}(X)$.  For each $ n $,
             for  $i \in \{0,1\}$  and  $ \vec{x} \in { \bf I }^{n} , T \in  {\cal S}_n(X) $  let
             us again use the linear maps   $e_{n,i}$,  \\
             \centerline{  $e_{n,i}: {\cal K}_{ n}(X)  \longrightarrow  {\cal K}_{ n}(X \times { \bf I })$,
                      $e_{n,i}(T) (\vec{x}) := ( T(\vec{x}), i )$, \  i.e.  $e_{n,i}(T) = e_{n,i} \circ T $, } \\
             see the definitions of $e_0$ and $e_1$ at the beginning of this section.  The boundary operator
             $ \partial_{n}$ is a  natural  transformation,  hence for  $i \in  \{0,1\}$
             the following diagram commutes, i.e.
             $  \partial_n \circ e_{n,i} =  e_{n-1,i} \circ \: \partial_n $ for all  $ n \in \mathbbm{N}$:
        \[   \xymatrix{   \cdots   \quad
            \ar[r]^<>(0.5){\partial_{n+2}}
            & {\cal K}_{ n+1}(X)     \ar[d]<0.0cm>_{e_{n+1,i}}   \ar[r]^<>(0.5){\partial_{n+1}}
          &  {\cal K}_{ n}(X)      \ar[d]<0.0cm>_{e_{n,i}}     \ar[r]^<>(0.5){\partial_{n}}
          &   {\cal K}_{ n-1}(X)   \ar[d]<0.0cm>_{e_{n-1,i}}   \ar[r]^<>(0.5){\partial_{n-1}}
          &   \quad     \cdots           \\
             \cdots  \quad
            \ar[r]_<>(0.5){\partial_{n+2}}
            &  {\cal K}_{ n+1}(X \times { \bf I })  \ar[r]_<>(0.5){\partial_{n+1}}
          &  {\cal K}_{ n}(X \times { \bf I })    \ar[r]_<>(0.5){\partial_{n}}
          &  {\cal K}_{ n-1}(X \times { \bf I })  \ar[r]_<(0.0)>(0.5){\partial_{n-1}}
          &                      \quad  \cdots        }
        \]
         The fact that $ \partial_{n}$ ({\sf u}) = 0   implies
             [$e_{n-1,i} \circ \,    \partial_{n}$] ({\sf u})  = 0 =
             [$   \partial_{n} \circ e_{n,i}$] ({\sf u}),
                and  we get that    $e_{n,i}$({\sf u})   is a cycle in   $X \times { \bf I }$,   thus                                [$e_{n,i}$({\sf u})]$_\sim \in  \, _{\vec{m}}{\cal H}_{n}(X \times { \bf I }$).

        Since we assumed that   $_{\vec{m}}{\cal H}$   satisfies the homotopy axiom,  by
                Lemma \ref{lemma eins}   the equivalence classes of   $ e_{n,0}$({\sf u}) and
                $ e_{n,1}$({\sf u}) in $  \, _{\vec{m}}{\cal H}_{n}(X \times { \bf I })$
                are the same.  This  means  we get the equality
                $[ e_{n,0}({\sf u}) ]_\sim = [ e_{n,1}({\sf u}) ]_\sim $,
                hence   $[ e_{n,0}({\sf u})  - e_{n,1}({\sf u})   ]_\sim = 0$.
                It follows that    $e_{n,0}( {\sf u} ) - e_{n,1}({\sf u})$  is  a  boundary,
                i.e.  there is a chain   $ \varphi \in  {\cal K}_{ n+1}(X \times { \bf I })$
                and   $e_{n,0}( {\sf u} ) - e_{n,1}( {\sf u} ) = \,    \partial_{n+1}( \varphi) . $
                Let $\varphi = \sum_{t=1}^{p}  r_t \cdot  \varphi_t  \in {\cal K}_{
n+1}(X \times { \bf I }) $ where $ p \in  \mathbbm{N}, \:   r_1 ,
r_2 ,  \ldots , r_p  \in {\cal R},                 \varphi_1 ,
\varphi_2 , \ldots ,   \varphi_p   \in  {\cal S}_{n+1}(X \times {
\bf I })$.
                Now  we  define four  maps
               $T_{1,0} , T_{2,0} , T_{1,1} , T_{2,1}  \in {\cal S}_{n}(X \times { \bf I })$. Let
         $$  T_{1,0} := e_{n,0} \circ T_{1} , \  T_{2,0} := e_{n,0} \circ T_{2} , \
                     T_{1,1} := e_{n,1} \circ T_{1} \ \text{and \ }   T_{2,1} := e_{n,1} \circ T_{2}. $$
              With    {\sf u}$  =  T_1 - T_2$    we have
              $e_{n,0}$({\sf u}) $ - \, e_{n,1}$({\sf u}) = $T_{1,0} - T_{2,0} - T_{1,1} + T_{2,1}$.
              Since   $T_1 \neq T_2$ the four maps $T_{1,0} , T_{2,0} , T_{1,1}$ and $T_{2,1} $                                         are  pairwise distinct.  We get:
    \begin{align*}
               T_{1,0} - T_{2,0} - T_{1,1} + T_{2,1}  =  e_{n,0}(T_{1} - T_{2}) -  e_{n,1}(T_{1} - T_{2})
               =   e_{n,0}( {\sf u} ) -  e_{n,1}( {\sf u} ) \\
               =   \partial_{n+1}( \varphi)
               =   \sum_{t=1}^{p}  r_t \cdot \,    \partial_{n+1}( \varphi_t)  =
               \sum_{t=1}^{p} r_t \cdot  \sum_{j =  1}^{n+1}   (-1)^{j+1} \cdot \sum_{i = 0}^{L}
               m_{i} \cdot \left\langle  \varphi_t\right\rangle_{ n+1, i, j } \, .
    \end{align*}
          The summands  $ \left\langle  \varphi_t\right\rangle_{ n+1, i, j } $
                 on the right hand side are elements of    ${\cal S}_{n}(X \times { \bf I })$
                (with coefficients $r_t, m_i$),   which generate the four summands
                 $ T_{1,0} , T_{2,0} , T_{1,1} , T_{2,1} $    on the left hand side.
                 Let us take the set  $B$  of triples,
       \begin{multline*}
      B :=\\ \{ (t , j , i ) \mid  t \in \{1,2, \ldots , p\} , j \in \{1,2, \ldots , n+1\}, i \in \{0,1, \ldots ,L \}
                 \wedge \left\langle  \varphi_t\right\rangle_{ n+1, i, j } =  T_{1,0} \}. \end{multline*}
                 Then we have
         \begin{align*}
              1_{\cal R}  \cdot T_{1,0}  =   \sum_{ ( t , j , i ) \in B }  r_t  \cdot  (-1)^{j+1}  \cdot
              m_{i}  \cdot \left\langle  \varphi_t\right\rangle_{ n+1, i, j }
                =   \sum_{ ( t , j , i ) \in B } r_t  \cdot  (-1)^{j+1}  \cdot   m_{i}  \cdot  T_{1,0} \ .
         \end{align*}
      This means that  $\Span _{\cal R}(\vec{m})  = {\cal R}$,
              and the proof of Theorem \ref{Homotopie Axiom} is complete.
   \end{proof}

   \section{The Exact Sequence of a Pair}
         As  we  mentioned  before, for all  $n  \in  \mathbbm{N}_0$  the boundary operator yields a functor
         $_{\vec{m}}{\cal H}_{n}:  \mathsf{TOP}^{2} \longrightarrow {\cal R}$-$ \mathsf{MOD}$,
         that means   for   any  $ ( f :  (X,A) \longrightarrow  (Y,B) )  \in   \mathsf{TOP}^{2}$
         we  have  a morphism of $ {\cal R}$-modules
         $_{\vec{m}}{\cal H}_{n}(f): \, _{\vec{m}}{\cal H}_{n}(X,A) \longrightarrow \, _{\vec{m}}{\cal H}_{n}(Y,B)$.

     For  any  cycle   {\sf u}  $ \in  {\cal K}_{  n}(X,A) $
         (i.e.   $   \partial_{n}(${\sf u})   is a chain in $A$, hence we have an equivalence class
         [{\sf u}]$_{\sim} \in \, _{\vec{m}}{\cal H}_{n}(X,A)$), we had abbreviated
         (at the end of section  \ref{section boundary operator})     \\
    \centerline{    $_{\vec{m}}{\cal H}_{n}(f) ($[{\sf u}]$_{\sim} )  =  f_\ast ( $[{\sf u}]$_{\sim} )   =
                          [{\cal K}_{  n}(f)$({\sf u})]$_{\sim} \in \, _{\vec{m}}{\cal H}_{n}(Y,B)$. }   \\
     For a  subspace $A \subset X$ we  get  a short exact  sequence of   ${\cal R}$-modules
    \[   \xymatrix{              \{ 0 \}                        \ar[r]
            &  {\cal K}_{ n}(A)                                   \ar[r]
          &  {\cal K}_{ n}(X)                                   \ar[r]
          &  {\cal K}_{ n}(X,A)                                 \ar[r]
          &                                       \{ 0 \}   .    }
    \]
   Together  with the  boundary operators    ($   \partial_{n } )_{n \geq 0}$  we  get a
          short exact sequence of  chain complexes
        \[   \xymatrix{              \{ 0 \}                            \ar[r]
            &  \, _{\vec{m}}{\cal K}_*(A)                               \ar[r]
          &  \, _{\vec{m}}{\cal K}_*(X)                               \ar[r]
          &  \, _{\vec{m}}{\cal K}_*(X,A)                             \ar[r]
          &                                       \{ 0 \}   .    }
        \]
      Let   $i: A \hookrightarrow X$    and    $j: (X, \emptyset ) \hookrightarrow (X, A)$  \
          be  the canonical topological inclusions. \
          Now we are able to construct for all $n \in \mathbbm{N}_0$ a  morphism  $k_\ast$
          of   ${\cal R}$-modules,  \\
     \centerline{     $k_\ast: \, _{\vec{m}}{\cal H}_{n}(X,A)  \longrightarrow \
                             _{\vec{m}}{\cal H}_{n-1}(A)$,      }    \\
          the {\it connecting homomorphism}. \
          Finally this  yields   a  long  exact  sequence   of  ${\cal R}$-module morphisms:
     \[ \cdots          \ \labto{ \  j_\ast  \ }  \,            _{\vec{m}}{\cal H}_{n+1}(X,A)                                                     \ \labto{ \  k_\ast  \ } \,            _{\vec{m}}{\cal H}_{n}(A)
                          \ \labto{ \  i_\ast  \ }  \,            _{\vec{m}}{\cal H}_{n}(X)
                        \ \labto{ \  j_\ast  \ }  \,            _{\vec{m}}{\cal H}_{n}(X,A)
                        \ \labto{ \  k_\ast  \ }  \             \cdots
     \]
           For   details   see   any   book  about singular homology theory,  for  instance
           \cite[p.93 ff]{Rotman}, or  \cite[p.18 ff]{Vick},  but  there is  no  necessity  for  us  to
           repeat  all  these   well known  facts.

\section{The Excision Axiom}
For the next section  it  is very useful  to compare the
corresponding section in  \cite[p.26 ff]{Massey}.  The reader should
note that we are able to prove the excision axiom only in the case
of   $ L=1$.  Furthermore, the reader may find perhaps an easier way
to prove the excision axiom, e.g. in \cite{Schoen}.    For a
topological space $X$ and a subset $A$,  $Int(A)$ means the interior
of $A$ and $Cl(A)$ the closure of $A$.  To prove  the following
theorem we need an extra assumption called ${\cal {NCD}}$, see
Definition \ref{definition drei}.
\begin{theorem}[Excision Axiom]  \label{Excision Axiom}
Let $X$ be a topological space and let $B$ and  $A$  be subsets of
$X$ such that $ Cl(B) \subseteq Int(A)$.  Let
 $$i: (X\backslash B , A\backslash B) \hookrightarrow (X, A)$$
be the canonical inclusion in $\mathsf{TOP}^{2}$.  Let us fix the
length $L := 1$, let the weight $\vec{m}$ be $\vec{m} := ( m_0, m_1
) := \left( a ,b \right) \in  {\cal R}^{2}$, let $(a,b)$ satisfy the
condition  ${\cal {NCD}}$. Then  for each $n  \in  \mathbbm{N}_0$
the morphism $$i_{\ast}: \, _{\vec{m}}{\cal H}_{n}( X \backslash B ,
A \backslash B ) \to \; _{\vec{m}}{\cal H}_{n}(X , A )$$  induced
by the inclusion $i$ is an isomorphism.
\end{theorem}

     This  theorem  follows directly from the following Proposition \ref{Proposition Isomorphismus}.
     Let  ${\cal U}$ := $\{ U_i  \mid i \in  \Im \}$
     be an indexed family of subsets of $X$ whose interiors cover $X$. Then $\cal U$ is called a
     {\it generalised open covering of  $X$. }  (The sets $U_i$ need  not be open).
     See  Definition \ref{definition vier}.      Note that the boundary operator $ \partial_{n}$
     commutes  with   the  inclusion  $\widehat{i}$, that means \,
     $ \partial_{n} \circ \widehat{i} (T) =   \widehat{i}  \circ \,    \partial_{n} (T)$
     for all   $T  \in    {\cal S}_{n}(A,\: {\cal U})$.  Therefore   $ \partial_{n}$  induces  a  linear  map
     ${\cal K}_{ n}(X, A, {\cal U}) \; \longrightarrow  {\cal K}_{ n-1 }(X, A, {\cal U}) $,
     which  we  call  $\partial_{n}$, too. Because $ \partial_{n} \circ \, \partial_{n+1} = 0$
     it leads  to  { \it  ${\cal U}$-small homology ${\cal R}$-modules }   \\
        \centerline{ $_{\vec{m}}{\cal H}_{n}(\: X , A , {\cal U} \:) :=  \frac {kernel ( \partial_{n })}
                     {image (\partial_{n+1})}$ . }   \\
     For   more  details see  \cite[p.29,30]{Massey}.  Now  we  are   able  to  formulate and to prove
     the  following  proposition.
  \begin{proposition}   \label{Proposition Isomorphismus}
       Let   $X$   be a topological  space  and  $A$ a subspace with the inclusion
                 $ A \stackrel{i}{\hookrightarrow}  X$,
                  and  let   ${\cal U}$   be a generalised open covering of $ X $ with the canonical inclusion
                  ${\cal K}_{ n}( X,A,{\cal U} )$  $\stackrel{j}{\hookrightarrow}$
                  ${\cal K}_{ n}( X,A )$.  Let us assume a weight  $\vec{m} := \left( a ,b \right) \in {\cal R}^{2}$
                  and let $ (a,b) $  satisfy the condition ${\cal {NCD}}$.  Then for all
                  $n  \in \mathbbm{N}_0$ the morphism
$$j_\ast \colon _{\vec{m}}{\cal H}_{n}( X,A, \: {\cal U})\stackrel{\cong }{\longrightarrow} _{\vec{m}}{\cal H}_{n}( X,A \: ) $$
induced by   $j$ is an  isomorphism.
   \end{proposition}
          The  proof  is rather  lengthy and  will  need  the  entire section.
          With  this   proposition the  excision axiom  easily  follows,  see  \cite[p.30,31]{Massey}.
          It  remains  to prove  the  proposition.
   \begin{proof} \quad
      First we shall present the proof with $A = \emptyset$, the empty set.
          Afterwards the general case  $ A \neq \emptyset$  is an easy application of the Five-Lemma.
          Hence let  $A := \emptyset $. \\
      We  have  to  define for all integers   $n \geq -1$   `subdivision maps'
$${\cal SD}_n  :    {\cal K}_{ n}(X)   \longrightarrow  {\cal K}_{ n}(X). $$

We need some preparations.   We define the ${\cal SD}_n$'s on the
basis $ {\cal S}_n(X) $  and  we  extend  the  definition on $ {\cal
K}_{ n}(X) $   by linearity.    For    $n = -1$, ${\cal SD}_{-1}$ is
the 0-map;  for $n = 0$  let    ${\cal SD}_0 (T) = -T $  for all $T
: { \bf I }^{0} \rightarrow  X $. For   $ n > 0$ a  map  $T  \in
{\cal S}_n(X) $    will be `subdivided'    into smaller ones.

Now we use the sets ${\cal E} := \{ 0 , 2 \} ,  {\cal V} :=  \{ -1,1
\}$. Further the reader should recall a map $\left\| T
\right\|_{\alpha , \vec{e} , \vec{v} }  \in  {\cal S}_n(X) $ of
Definition \ref{definition fuenf}. In the following we set  $\alpha
:= \frac{1}{3}$ .  We take  $ e_i  \in {\cal E}$  and $ v_i \in
{\cal V}$ for $i= 1,2, \ldots , n $,  and $ q_n $  will be the
identity on  ${ \bf I }^{n}$. Define  for all  $n \in \mathbbm{N}$
and  for any continuous $T: { \bf I }^{n} \rightarrow  X \
(\text{i.e. }  T \in  {\cal S}_{n}(X) ):$
\begin{align}     \label{definition SD}
{\cal SD}_n  (T)   :=        \sum_{\vec{e} \in {\cal E}^{n}  }   {}\
\sum_{\vec{v} \in {\cal V}_{\vec{e},n}} {}\ ( - \prod_{i=1}^{n}  v_i
) \cdot \left\| T \right\|_{\frac{1}{3}, \vec{e}, \vec{v}}  \: ,
\end{align}
where $ {\cal V}_{\vec{e},n}$ is the set of $( v_1, v_2,  \ldots  ,
v_n ) \in {\cal V}^{n}$ such that for all   $ i = 1, 2, \ldots , n
$, $v_i = 1$ if $ e_i = 0 $ and $ v_i \in \{-1,1 \}$  if    $ e_i =
2 $.      We  get a map  ${\cal SD}_n : {\cal K}_{ n}(X)
\longrightarrow {\cal K}_{ n}(X) $  by linearity.

\noindent Examples:  Let $ n := 1 $. For $T \in  {\cal S}_{1}(X)$,
i.e. $T: { \bf I} \rightarrow X $ we have
$${\cal SD}_1  (T) \ =  \ -  \left\| T \right\|_{\frac{1}{3} , 0 , 1 } +
\left\| T \right\|_{\frac{1}{3} , 2 , -1 } -  \left\| T
\right\|_{\frac{1}{3} , 2 , 1 } \quad \text{(see Figure \ref{pic3})}
,
$$ and  for   $ n := 2 $, for $T  \in {\cal S}_{2}(X)$   we get the
linear combination (see  the figure, too)
   \begin{multline*}
{\cal SD}_2  (T)  = -  \left\| T \right\|_{\frac{1}{3} , {0 \choose
0} , {1 \choose 1} } -  \left\| T \right\|_{\frac{1}{3} , {2 \choose
0} , {1 \choose 1} } +  \left\| T \right\|_{\frac{1}{3} , {2 \choose
0} , {-1 \choose 1} } -  \left\| T \right\|_{\frac{1}{3} , {0 \choose 2} , {1 \choose 1} }    \\
+  \left\| T \right\|_{\frac{1}{3} , {0 \choose 2} , {1 \choose -1}
}  +  \left\| T \right\|_{\frac{1}{3} , {2 \choose 2} , {-1 \choose
1} }  -  \left\| T \right\|_{\frac{1}{3} , {2 \choose 2} , { 1
\choose 1} }  +  \left\| T \right\|_{\frac{1}{3} , {2 \choose 2} ,
{1 \choose -1} } -  \left\| T \right\|_{\frac{1}{3} , {2 \choose 2}
, {-1 \choose -1} } \, _{.}
\end{multline*}
   \begin{figure}[ht]
        \centering
        \setlength{\unitlength}{1cm}
        \begin{picture}(4,4.5)
          \put(-1.5,2.9 ){ $ {\cal SD}_1  (T):$ }
          \put( 4.75,3.4 ){ $ {\cal SD}_2  (T):$ }
          \put(-2,1){\vector(1,0){0.5}}     \put(-2,0.9){\line(0,1){0.2}}                             \put(-2.1,0.55){0}
          \put(1,1){\vector(-1,0){1.5}}     \put(-1,0.9){\line(0,1){0.2}} 
          \put(-1.1,0.5 ){$\frac{1}{3}$}
          \put(-2,1){\vector(1,0){2.5}}     \put(0,0.9) {\line(0,1){0.2}} 
          \put(-0.1,0.5 ){$\frac{2}{3}$}    \put(1,0.9) {\line(0,1){0.2}}  
                                            \put(0.9,0.5 ){1}

          \put(-4.2,2.0 ){ $-  \left\| T \right\|_{\frac{1}{3} , 0 , 1 }$ }   \put(-2.8,1.7){\line(2,-1){1.1}}
          \put(-2.0,2.0 ){ $ + \left\| T \right\|_{\frac{1}{3} , 2 , -1 }$ }  \put(-0.8,1.7){\line(1,-2){0.3}}
          \put(0.2,2.0 ){ $ -  \left\| T \right\|_{\frac{1}{3} , 2 , 1 }$ }   \put(0.4,1.2){\line(1,1){0.6}}
          \put(4,0){\vector(1,0){0.5}}      \put(4.35,0.4){$-$}      \put(3.9,-0.4){0}
          \put(4,0){\vector(1,0){2.5}}      \put(5.35,0.4){$+$}      \put(4.9,-0.4){$\frac{1}{3}$}
          \put(7,0){\vector(-1,0){1.5}}     \put(6.35,0.4){$-$}      \put(5.9,-0.4){$\frac{2}{3}$}
          \put(4,1){\vector(1,0){0.5}}                               \put(6.9,-0.4){$1$}
          \put(4,1){\vector(1,0){2.5}}
          \put(7,1){\vector(-1,0){1.5}}
                                  \put(4.35,1.4){$+$}         \put(3.7,0.9){$\frac{1}{3}$}
                                  \put(5.35,1.4){$-$}         \put(3.7,1.9){$\frac{2}{3}$}
                                  \put(6.35,1.4){$+$}         \put(3.7,2.9){$1$}
           \put(4,2){\vector(1,0){0.5}}
           \put(4,2){\vector(1,0){2.5}}
           \put(7,2){\vector(-1,0){1.5}}
           \put(4,3){\vector(1,0){0.5}}
           \put(4,3){\vector(1,0){2.5}}
           \put(7,3){\vector(-1,0){1.5}}
           \put(4,0){\vector(0,1){0.5}}
           \put(4,0){\vector(0,1){2.5}}
           \put(4,3){\vector(0,-1){1.5}}
                                 \put(4.35,2.4){$-$}
                                 \put(5.35,2.4){$+$}
                                 \put(6.35,2.4){$-$}
           \put(5,0){\vector(0,1){0.5}}
           \put(5,0){\vector(0,1){2.5}}
           \put(5,3){\vector(0,-1){1.5}}
           \put(6,0){\vector(0,1){0.5}}
           \put(6,0){\vector(0,1){2.5}}
           \put(6,3){\vector(0,-1){1.5}}
           \put(7,0){\vector(0,1){0.5}}
           \put(7,0){\vector(0,1){2.5}}
           \put(7,3){\vector(0,-1){1.5}}
      \end{picture}
\caption{} \label{pic3}
  \end{figure}

Generally   for    $n  \in \mathbbm{N}_0$   and    $T \in   {\cal
S}_{n}(X)$, ${\cal SD}_n (T) $    is a linear combination  of
$3^{n}$   maps in $ {\cal S}_{n}(X) $.
\begin{lemma}  \label{lemma drei}   \quad
For all   $n  \in \mathbbm{N}_0$   the  map  ${\cal SD}_n$  commutes
with the  boundary operator   $   \partial_{n}$,   i.e.  the
following  diagram commutes:
\[   \xymatrix{   \cdots \quad    \ar[r]^<>(0.4){\partial_{n+2}}
            & {\cal K}_{ n+1}(X)  \ar[d]<0.0cm>^<>(0.5){ {\cal SD}_{n+1}}
                                        \ar[r]^<>(0.5){\partial_{n+1}}
          &  {\cal K}_{ n}(X)   \ar[d]<0.0cm>^<>(0.5){ {\cal SD}_{n}}
                                      \ar[r]^<>(0.5){\partial_{n}}
          &  {\cal K}_{ n-1}(X) \ar[d]<0.0cm>^<>(0.5){ {\cal SD}_{n-1}}
                                      \ar[r]^<(0.4)>(0.5){\partial_{n-1}}
          &  {\cal K}_{ n-2}(X) \ar[d]<0.0cm>^<>(0.5){ {\cal SD}_{n-2}}
                                      \ar[r]^<>(0.5){\partial_{n-2}}
          &                        \quad        \cdots           \\
              \cdots   \quad  \ar[r]_<>(0.4){\partial_{n+2}}
            &  {\cal K}_{ n+1}(X)   \ar[r]_<>(0.5){\partial_{n+1}}
          &  {\cal K}_{ n}(X)     \ar[r]_<>(0.5){\partial_{n}}
          &  {\cal K}_{ n-1}(X)   \ar[r]_<(0.0)>(0.5){\partial_{n-1}}
          &  {\cal K}_{ n-2}(X)   \ar[r]_<>(0.5){\partial_{n-2}}
          &                                            \quad   \cdots      \\     }
        \]
\end{lemma}
\begin{proof} \quad
We  have  to prove $ \partial_{n} \circ  {\cal SD}_{n} (T) = {\cal
SD}_{n-1} \circ \, \partial_{n} (T)$ for each $n \in \mathbbm{N}_0 $
and $T \in  {\cal S}_{n}(X) $. This is trivial for $n = 0$  and easy
for $n = 1$, so let $n  \geq 2$.  Let $T \in {\cal S}_{n}(X) $.  \
Note that in the following  we shall use    `$\left\langle
T\right\rangle_{\: i,\; j \;} $' instead of the expression
`$\left\langle T\right\rangle_{\: n,\;  i,\; j \;} $', to make it
better readable.  We have:
    \begin{align*}
         &  \partial_{n} \circ  {\cal SD}_{n} (T)  =  \partial_{n}  \left[ \sum_{\vec{e} \in {\cal E}^{n} } \
                    \sum_{\vec{v} \in {\cal V}_{\vec{e},n}} {} ( - \prod_{i=1}^{n}  v_i ) \cdot
                    \left\| T \right\|_{\frac{1}{3} , \vec{e} , \vec{v} } \quad \right]          \\
         &  =   \sum_{j=1}^{n} (-1)^{j+1}  \sum_{\vec{e} \in {\cal E}^{n}  } \
                    \sum_{\vec{v} \in {\cal V}_{\vec{e},n}} {}\ ( - \prod_{i=1}^{n} v_i ) \cdot \left[   a \cdot
                    \left\langle \left\| T \right\|_{\frac{1}{3} , \vec{e} , \vec{v} } \right\rangle_{0,j}   + b \cdot
                    \left\langle \left\| T \right\|_{\frac{1}{3} , \vec{e} , \vec{v} } \right\rangle_{1,j} \right]  .
    \end{align*}
   As well as
   \begin{align*}
         &   {\cal SD}_{n-1} \circ \partial_{n} (T)  =  {\cal SD}_{n-1}  \left[ \sum_{j=1}^{n} (-1)^{j+1}
                    \left( a \cdot \left\langle T \right\rangle_{0,j} + b \cdot
                    \left\langle T \right\rangle_{1,j} \right) \quad \right]                \\
         &  =      \sum_{j=1}^{n} (-1)^{j+1}  \sum_{\vec{e} \in {\cal E}^{n-1}  } {}\
                    \sum_{\vec{v} \in {\cal V}_{\vec{e},n-1}}  {}
                    (- \prod_{i=1}^{n-1}  v_i) \cdot \left[    a \cdot  \left\|  \left\langle T \right\rangle_{0,j}
                    \right\|_{\frac{1}{3} , \vec{e} , \vec{v} }    +
                    b \cdot  \left\|  \left\langle T \right\rangle_{1,j}
                    \right\|_{\frac{1}{3} , \vec{e} , \vec{v} }    \right]  .
   \end{align*}
     The equality is not obvious; so we have to calculate. It seems that the first sum is `bigger'. But many
           elements cancel pairwise, and the rest is equal to  the second sum.

     We fix an arbitrary $j \in \{ 1 , 2 , \ldots , n\}$.  We take
$$\vec{e_2} =  ( e_1,\: e_2, \ldots , e_{j-1},\: 2 ,\: e_{j+1},  \ldots , e_n
  )\in  {\cal E}^{n}$$
and $\vec{\vartheta_1}, \vec{\vartheta_{-1}} \in {\cal
V}_{\vec{e_2},n}$ given by
$$\vec{\vartheta_1} = ( v_1, v_2, \ldots , v_{j-1}, \: 1 ,\: v_{j+1},  \ldots , v_n ),
\quad \vec{\vartheta_{-1}} = ( v_1, v_2, \ldots , v_{j-1},\: -1 ,\:
v_{j+1}, \ldots  , v_n ) .$$ Then $$\left\langle \left\| T
\right\|_{\frac{1}{3} , \vec{e_2} , \vec{\vartheta_1} }
\right\rangle_{0,j} ,  \left\langle \left\| T \right\|_{\frac{1}{3}
, \vec{e_2} ,
 \vec{\vartheta_{-1}} } \right\rangle_{0,j}                   \in  {\cal S}_{n-1}(X) . $$
     For a point   $( x_1, x_2, \ldots , x_{j-1}, x_{j}, \ldots , x_{n-1})  \in { \bf I }^{n-1}$  we get
     \begin{align*}
           &   \left\langle \left\| T \right\|_{\frac{1}{3} , \vec{e_2} , \vec{\vartheta_1} } \right\rangle_{0,j}
           ( x_1, x_2, \ldots , x_{j-1} , x_{j}, \ldots ,  x_{n-1} )  \\
        =  & \left\| T \right\|_{\frac{1}{3} , \vec{e_2} , \vec{\vartheta_1} }
           ( x_1, x_2, \ldots , x_{j-1} ,\: 0 ,\: x_{j} , \ldots,  x_{n-1}) \\
        =  & T \left( \frac{1}{3} \cdot [e_1+v_1\cdot x_1, \ldots ,e_{j-1}+v_{j-1}\cdot x_{j-1}, \: 2 , \:
            e_{j+1}+v_{j+1}\cdot x_{j} ,  \ldots , e_{n}+v_{n}\cdot x_{n-1} ] \right) ,
     \end{align*}
and \begin{align*}  \left\langle \left\| T \right\|_{\frac{1}{3} ,
\vec{e_2} , \vec{\vartheta_{-1}} } \right\rangle_{0,j} &=
  \left\| T \right\|_{\frac{1}{3} , \vec{e_2} ,
\vec{\vartheta_{-1}} }
            ( x_1,\: \ldots ,\: x_{j-1} \:,\: 0\: ,\: x_{j}\: ,\: \ldots \:,\:
            x_{n-1})\\
&=  T \left( \frac{1}{3} \cdot [e_1+v_1\cdot x_1, \: \ldots  \ldots
\;, 2 \: , \
              \: \ldots  \ldots \;, e_{n}+v_{n}\cdot x_{n-1}\: ]
              \right)\\
\intertext{whence} \left\langle \left\| T \right\|_{\frac{1}{3} ,
\vec{e_2} , \vec{\vartheta_1} } \right\rangle_{0,j}
           &= \left\langle \left\| T \right\|_{\frac{1}{3} , \vec{e_2}, \vec{\vartheta_{-1}}} \right\rangle_{0,j}.
 \end{align*}
 Now note that \ $ ( \prod_{v_i \in \vec{\vartheta_1}} v_i ) \cdot ( \prod_{v_i \in \vec{\vartheta_{-1}}} v_i )
              = -1 $,   from which  it follows that
\begin{align*}
a \cdot (   -\prod_{v_i \in \vec{\vartheta_1}} \:  v_i   ) \cdot  \
\left\langle \left\| T \right\|_{\frac{1}{3} , \vec{e_2} ,
\vec{\vartheta_1} } \right\rangle_{0,j} +   a \cdot (   -\prod_{v_i
\in \vec{\vartheta_{-1}}} \:  v_i )  \cdot  \ \left\langle \left\| T
\right\|_{\frac{1}{3} , \vec{e_2} , \vec{\vartheta_{-1}} }
\right\rangle_{0,j} \ = \  0   .
\end{align*}
In the same way  with
\begin{alignat*}{2}
\vec{e_0} & := ( e_1, e_2, \ldots , e_{j-1},\: 0 ,\: e_{j+1}, \ldots
, e_n )   ,& \vec{e_2} & := ( e_1, e_2, \ldots , e_{j-1},\: 2 ,\:
e_{j+1}, \ldots , e_n )\in  {\cal E}^{n}\\
\vec{\vartheta_1} & := ( v_1,  \ldots , v_{j-1}, \: 1 ,\: v_{j+1},
\ldots , v_n ),& \vec{\vartheta_{-1}} &:= ( v_1, \ldots , v_{j-1},\:
-1 ,\: v_{j+1}, \ldots , v_n ) \in  {\cal V}_{\vec{e_2},n}
\end{alignat*}
we get
\begin{align*}
&\left\langle \left\| T \right\|_{\frac{1}{3} , \vec{e_0} ,
\vec{\vartheta_1} } \right\rangle_{1,j} ( x_1, x_2, \ldots , x_{n-1}
)  \\ &=   \left\| T \right\|_{\frac{1}{3}, \vec{e_0},
\vec{\vartheta_1}} ( x_1, x_2, \ldots , x_{j-1} ,\: 1 ,\: x_{j} ,
\ldots, x_{n-1})\\
&=T \left( \frac{1}{3} \cdot [e_1+v_1\cdot x_1,  \ldots ,
e_{j-1}+v_{j-1} \cdot x_{j-1}, \: 1 , \: e_{j+1}+v_{j+1} \cdot
x_{j}, \ldots , e_{n}+v_{n}\cdot x_{n-1}] \right)\\
&= \left\| T \right\|_{\frac{1}{3} , \vec{e_2} ,
\vec{\vartheta_{-1}} }( x_1, x_2, \ldots , x_{j-1} ,\: 1 ,\: x_{j} ,
\ldots, x_{n-1})\\
&= \left\langle \left\| T \right\|_{\frac{1}{3} , \vec{e_2} ,
\vec{\vartheta_{-1}} } \right\rangle_{1,j}  ( x_1, x_2, \ldots ,
x_{n-1} ) .
\end{align*}
Hence
\begin{align*}
b \cdot (   -\prod_{v_i \in \vec{\vartheta_1}} \:  v_i   )  \cdot  \
\left\langle \left\| T \right\|_{\frac{1}{3} , \vec{e_0} ,
\vec{\vartheta_1} } \right\rangle_{1,j} +   b \cdot (   -\prod_{v_i
\in \vec{\vartheta_{-1}}} \:  v_i ) \cdot \ \left\langle \left\| T
\right\|_{\frac{1}{3} , \vec{e_2} , \vec{\vartheta_{-1}} }
\right\rangle_{1,j}   &=  0   .
\end{align*}

Now take again
$$\vec{e_0} = ( e_1, \ldots , e_{j-1},\: 0 ,\: e_{j+1}, \ldots , e_n
) \in  {\cal E}^{n}, \quad \vec{\vartheta_1} = ( v_1, \ldots ,
v_{j-1},\: 1 ,\: v_{j+1}, \ldots  , v_n ) \in  {\cal
V}_{\vec{e_0},n}$$
 as above, and  define
 $$\widetilde{e} = ( e_1, e_2, \ldots , e_{j-1}, e_{j+1}, \ldots , e_n )
 \in  {\cal E}^{n-1}, \quad  \widetilde{\vartheta}   =
 ( v_1, v_2, \ldots , v_{j-1}, v_{j+1},  \ldots  , v_n ) \in {\cal V}_{\widetilde{e},n-1}.$$
Then   $$\left\langle \left\| T \right\|_{\frac{1}{3} , \vec{e_0} ,
\vec{\vartheta_1} }      \right\rangle_{0,j} ,      \left\|
\left\langle  T  \right\rangle_{0,j} \right\|_{\frac{1}{3},
\widetilde{e}, \widetilde{\vartheta}} \ \in  {\cal S}_{n-1}(X) ,$$
and  for all  points    $( x_1, x_2, \ldots , x_{j-1}, x_{j}, \ldots
, x_{n-1} )  \in { \bf I }^{n-1}$ we calculate
\begin{align*}
  &\left\langle \left\| T \right\|_{\frac{1}{3} , \vec{e_0} ,
\vec{\vartheta_1} } \right\rangle_{0,j} ( x_1, x_2, \ldots , x_{n-1}
)\\ &=  \left\| T \right\|_{\frac{1}{3} ,
          \vec{e_0} , \vec{\vartheta_1}} ( x_1, x_2, \ldots , x_{j-1}, \: 0 ,\: x_{j} , \ldots , x_{n-1})    \\
& =    T \left( \frac{1}{3} \cdot [e_1+v_1\cdot x_1,  \ldots ,
e_{j-1}+v_{j-1}\cdot x_{j-1}, \: 0 , \:e_{j+1}+v_{j+1}\cdot x_{j} , \ldots , e_{n}+v_{n}\cdot x_{n-1} ]   \right)   \\
& =   \left\langle  T  \right\rangle_{0,j} \left( \frac{1}{3}
\cdot [e_1+v_1\cdot x_1, \ldots , e_{j-1}+v_{j-1}\cdot x_{j-1}, \,
e_{j+1}+v_{j+1}\cdot x_{j} ,  \ldots , e_{n}+v_{n}\cdot x_{n-1} ] \right)  \\
& =  \left\| \left\langle  T  \right\rangle_{0,j}
\right\|_{\frac{1}{3}, \widetilde{e}, \widetilde{\vartheta}}
( x_1, x_2,  \ldots , x_{j-1} , x_{j} , \ldots , x_{n-1} )  .
     \end{align*}
Hence
\begin{align*}
            a \cdot (   -\prod_{v_i \in \vec{\vartheta_1}} \:  v_i ) \cdot \
            \left\langle \left\| T \right\|_{\frac{1}{3} , \vec{e_0} , \vec{\vartheta_1} } \right\rangle_{0,j}   = \
            a \cdot (   -\prod_{v_i \in \widetilde{\vartheta}} \:  v_i ) \cdot \ \left\| \left\langle  T
             \right\rangle_{0,j} \right\|_{\frac{1}{3}, \widetilde{e}, \widetilde{\vartheta}}  \ .
\end{align*}

If we take   as above,
           $$\vec{e_2} = ( e_1, e_2, \ldots , e_{j-1},\: 2 ,\: e_{j+1},  \ldots , e_n ) , \
           \widetilde{e} = ( e_1, e_2, \ldots , e_{j-1}, e_{j+1}, \ldots , e_n )$$
           and
$$\vec{\vartheta_1} = (v_1, v_2, \ldots , v_{j-1}, \: 1 ,\: v_{j+1}, \ldots , v_n ) \in {\cal
V}_{\vec{e_2},n}\quad  \widetilde{\vartheta} = ( v_1, v_2, \ldots ,
v_{j-1}, v_{j+1}, \ldots , v_n )  \in
           {\cal V}_{\widetilde{e},n-1}$$
we have
$$\left\langle \left\| T \right\|_{\frac{1}{3} , \vec{e_2} , \vec{\vartheta_1} } \right\rangle_{1,j}, \
\left\| \left\langle  T  \right\rangle_{1,j} \right\|_{\frac{1}{3},
\widetilde{e}, \widetilde{\vartheta}} \in  {\cal S}_{n-1}(X) .$$
 We compute
\begin{align*}
&\left\langle \left\| T \right\|_{\frac{1}{3} , \vec{e_2} ,
\vec{\vartheta_1} } \right\rangle_{1,j}( x_1, x_2, \ldots ,
x_{n-1} ) \\
& =  \left\| T \right\|_{\frac{1}{3} , \vec{e_2} ,
\vec{\vartheta_1}}             ( x_1, x_2, \ldots , x_{j-1},\: 1 ,\: x_{j} ,  \ldots , x_{n-1} )    \\
& =   T \left( \frac{1}{3}   \cdot [e_1+v_1\cdot x_1, \ldots ,
e_{j-1}+v_{j-1}\cdot x_{j-1},\: 3 , \: e_{j+1}+v_{j+1}\cdot x_{j}, \ldots , e_{n}+v_{n}\cdot x_{n-1} ] \right) \\
& =  \left\langle  T  \right\rangle_{1,j} \left( \frac{1}{3}   \cdot
[e_1+v_1\cdot x_1 , \ldots , e_{j-1}+v_{j-1}\cdot x_{j-1},
e_{j+1}+v_{j+1}\cdot x_{j} , \ldots , e_{n}+v_{n}\cdot x_{n-1} ]\right)\\
&= \left\| \left\langle  T  \right\rangle_{1,j}
\right\|_{\frac{1}{3}, \widetilde{e}, \widetilde{\vartheta}}( x_1,
x_2, \ldots , x_{j-1}, x_{j}, \ldots , x_{n-1} ).
\end{align*}
 Hence  we get
\begin{align*}     \quad
b \cdot (   -\prod_{v_i \in \vec{\vartheta_1}} \:  v_i ) \cdot \
\left\langle \left\| T \right\|_{\frac{1}{3} , \vec{e_2} ,
\vec{\vartheta_1} } \right\rangle_{1,j}  & = \ b \cdot ( -\prod_{v_i
\in \widetilde{\vartheta}} \:  v_i ) \cdot \ \left\| \left\langle  T
\right\rangle_{1,j} \right\|_{\frac{1}{3}, \widetilde{e},
\widetilde{\vartheta}} \,  .
\end{align*}
    This is all  we need to show that $ \partial_{n} \circ  {\cal SD}_{n} (T) =
             {\cal SD}_{n-1} \circ  \partial_{n} (T)$, and Lemma  \ref{lemma drei} has been proved.
  \end{proof}
    Because of  Lemma  \ref{lemma drei}, the map ${\cal SD}_n : {\cal K}_{ n}(X) \longrightarrow {\cal K}_{ n}(X) $
         induces an ${\cal R}$-module endomorphism  of the space $_{\vec{m}}{\cal H}_{n}(X)$, which we
         also call   ${\cal SD}_n$.

    Now  we shall show that for the weight $\vec{m} =  (a, b) $  we have
          $$ [ a \cdot {\cal SD}_n ({\sf u}) ]_{\sim} =  [ b \cdot {\sf u}   ]_{\sim} $$
          on the level of homology classes, for all
          {\sf u} $\in   {\cal K}_{   n}(X)$  with  {\sf u} $\in \text{ kernel} (\partial_{n})$.
     We are able to do this by the help of a chain homotopy in the same way we used it for the proof of
          the  homotopy axiom.    This means  for   $n \geq -1$    the construction  of
          a  linear   map  $\Theta_n :  {\cal K}_{ n}(X) \longrightarrow {\cal K}_{ n+1}(X)$,
          which yields  a  (noncommutative) diagram
     \[   \xymatrix{   \cdots \qquad
            \ar[r]^<>(0.4){\partial_{n+2}}
            & {\cal K}_{ n+1}(X)  \ar[dd]<0.1cm>^<>(0.6){Id}
                                        \ar[dd]<-0.15cm>_<>(0.6){{\cal SD}_{n+1}}
                                      \ar@<0.0 cm>[ldd]_<>(0.4){ \Theta_{n+1}}  \ar[r]^<>(0.5){\partial_{n+1}}
          &  {\cal K}_{ n}(X)   \ar[dd]<0.1cm>^<>(0.6){Id}
                                      \ar[dd]<-0.15cm>_<>(0.6){{\cal SD}_{n}}
                                      \ar[ldd]_<>(0.4){ \Theta_{n}} \ar[r]^<>(0.5){\partial_{n}}
          &  {\cal K}_{ n-1}(X) \ar[dd]<0.1cm>^<>(0.6){Id}
                                      \ar[dd]<-0.15cm>_<>(0.6){{\cal SD}_{n-1}}
                                      \ar[ldd]_<>(0.4){ \Theta_{n-1}}  \ar[r]^<(0.4)>(0.5){\partial_{n-1}}
          &  {\cal K}_{ n-2}(X) \ar[dd]<0.1cm>^<>(0.6){Id}
                                      \ar[dd]<-0.15cm>_<>(0.6){{\cal SD}_{n-2}}
                                      \ar[ldd]_<>(0.4){ \Theta_{n-2}}  \ar[r]^<>(0.5){\partial_{n-2}}
          &                     \ar[ldd]_<>(0.4){ \Theta_{n-3}}     \quad        \cdots           \\
                      \cdots  \qquad     & & & & &    \quad \cdots                                    \\
          \cdots \qquad
            \ar[r]_<>(0.4){\partial_{n+2}}
            &  {\cal K}_{ n+1}(X)  \ar[r]_<>(0.5){\partial_{n+1}}
          &  {\cal K}_{ n}(X)  \ar[r]_<>(0.5){\partial_{n}}
          &  {\cal K}_{ n-1}(X) \ar[r]_<(0.0)>(0.5){\partial_{n-1}}
          &  {\cal K}_{ n-2}(X)   \ar[r]_<>(0.5){\partial_{n-2}}
          &                                                                    \quad   \cdots      \\     }
        \]    \\
      Let  $Id $  be the identity map on  ${\cal K}_{ n}(X)$, for all $ n $.  Our aim is,  for
      ${\sf u} \in   {\cal K}_{  n}(X)$,   to get the equation
  \begin{align}   \label{equation sieben}
      ( \partial_{n+1} \circ \Theta_n )( {\sf u} )  =  \pm( b \cdot Id - a \cdot {\cal SD}_n ) ( {\sf u} )
      +  ( \Theta_{n-1} \circ  \partial_{n} ) ( {\sf u} ),  \ { \rm for } \ n \in \mathbbm{N}_{0} .
  \end{align}
Of course   $\Theta_{-1} := 0$. For  $ n := 0$  for every  $T: \{0\}
\rightarrow X$  and  for $x \in  { \bf I } = [0,1] $   define
$\Theta_0(T) (x) := T(0)$.   Then   ($  \partial_{1} \circ \Theta_0
)(T)  =  \ a \cdot T + b \cdot T  = +( b \cdot T - a \cdot {\cal
SD}_0(T) $),  as required. Let $n \geq 1$. We need three auxiliary
functions $\eta_0 , \eta_1 , \eta_2 : \ {\bf I}^{2} \rightarrow {
\bf I }$: for all  $x,y \in [ 0,1 ]$ let
\begin{align*}
\eta_0 (x,y) &:= \frac{x}{3-2 \cdot y},\\
 \eta_1(x,y)  &  :=     \begin{cases}  \frac{2 - x}{3 - 2 \cdot y}    & \text{  for }
 y  \leq  \frac{1}{2} +  \frac{1}{2} \cdot x  ,    \\
 1 & \text{  otherwise }    ,
 \end{cases}\\
\eta_2(x,y)   & :=     \begin{cases}  \frac{2 + x}{3 - 2 \cdot y} &
\text{  for }
y  \leq  \frac{1}{2} -  \frac{1}{2} \cdot x  ,    \\
1 & \text{  otherwise }   .
\end{cases}
\end{align*}
The maps   $\eta_0 , \eta_1 , \eta_2 $   are continuous. We use the
set  $\Upsilon := \{ 0, 1, 2 \}$.
     For all tuples \: $\vec{z} = ( z_1, z_2, \ldots , z_n ) \in \Upsilon^{n}$ and all $T \in {\cal S}_{n}(X) $
               define \:  $G_{\vec{z}}(T)  : { \bf I }^{n+1} \rightarrow X$  by the equation
      $$G_{\vec{z}}(T) ( x_1, \ldots , x_n, x_{n+1} ) := \
            T ( \eta_{z_1}(x_1,x_{n+1}), \eta_{z_2}(x_2,x_{n+1}), \ldots \ldots , \eta_{z_n}(x_n,x_{n+1}) ) \, , $$
      for all  $(n+1)$-tuples   $( x_1,  \ldots ,   x_n, x_{n+1} ) \in { \bf I }^{n+1}$.
      $G_{\vec{z}}(T)$   is  an element of $ {\cal S}_{n+1}(X) $.
      Let for all   $\vec{z} := ( z_1, z_2, \ldots , z_n ) \in \Upsilon^{n}$ the number
      $v_{\vec{z}}  \in \{ -1, +1 \} $   by
   $v_{\vec{z}}  := (-1)^{\sum_{i=1}^{n} z_i } $, hence
      $v_{\vec{z}} = (-1)^\alpha$ where $\alpha$ is the  number  of  1's  in  $\vec{z} $,  and finally define:
   \begin{align}
           \Theta_n (T) := \sum_{ \vec{z} \in \Upsilon^{n} } \quad  v_{\vec{z}} \cdot  G_{\vec{z}}(T)  .
   \end{align}
    We describe an example for $n = 1$. For  $T \in {\cal S}_{1}(X) $ we have for all pairs
          $(x,y) \in { \bf I }^{2}$:
    \begin{align*}
              \Theta_1(T) (x,y)  &  = + G_0(T)(x,y) - G_1(T)(x,y) + G_2(T)(x,y)      \\
               &  = + T( \eta_0(x,y) ) - T( \eta_1(x,y) ) + T( \eta_2(x,y) ) \, .
    \end{align*}
   We hope that the following Figure \ref{pic4} will help to get a better understanding. By definition, $\Theta_1(T)$
      generates three maps   $ T(\eta_0), T(\eta_1), T(\eta_2) \in  {\cal S}_{2}(X) $, whose images are indicated
      by squares.  Later we shall prove that the top of the left square equals the image of  $T$,   while
      the bottoms of all squares are equal to    $-\;{\cal SD}_1 (T)$. Further, the upper fourth (diagonal) section of
      the  middle square is constant equal $-T(1)$, and the upper three fourths section of the right square is constant
      equal   $T(1)$.
       \\  \\  \\  \\
     \begin{figure}[h]
          \centering
          \setlength{\unitlength}{1cm}
             \begin{picture}(6,4.5)
                    \put(-3,3.5){$\Theta_1(T)$ = }
                    \put(-1,2){\line(0,1){3}}     \put(5.5,2){\line(0,1){3}}
                    \put(2,2){\line(0,1){3}}      \put(6,2){\line(0,1){3}}
                    \put(2.5,2){\line(0,1){3}}    \put(9,2){\line(0,1){3}}
                    \put(-1,2){\line(1,0){3}}     \put(-1,5){\line(1,0){3}}
                    \put(2.5,3.5){\line(2,1){3}}  \put(6,3.5){\line(2,-1){3}}
                    \put(1.1,5.3){$ T$ }          \put(0.48,5.1){\line(2,1){0.5}}
                    \put(2.5,0.2){$-\;{\cal SD}_1 (T)$ }

                    \put(3.2,0.6){\line(-2,1){2.5}}   \put(3.4,0.6){\line(1,4){0.32}}   \put(3.6,0.6){\line(3,1){3.9}}
                    \put(-0.2,2.7){$T(\eta_0(x,y))$}  \put(3.0,2.7){$-T(\eta_1(x,y))$}  \put(6.7,2.6){$T(\eta_2(x,y))$}
                    \put(3.1,4.5){$ -T(1)$}     \put(7.25,3.9){$ T(1)$}
                    \thicklines \put(-1,4.95){\line(1,0){3}}  \put(-1,4.96){\line(1,0){3}}  \put(-1,4.97){\line(1,0){3}}
                    \put(-1,4.98){\line(1,0){3}}        \put(-1,4.99){\line(1,0){3}}
                    \put(-1,2.04){\line(1,0){3}}        \put(-1,2.03){\line(1,0){3}}
                    \put(-1,2.02){\line(1,0){3}}        \put(-1,2.01){\line(1,0){3}}
                    \put(2.5,2.04){\line(1,0){3}}      \put(2.5,2.03){\line(1,0){3}}     \put(2.5,2.02){\line(1,0){3}}                      \put(2.5,2.01){\line(1,0){3}}
                    \put(6,2.04){\line(1,0){3}}        \put(6,2.03){\line(1,0){3}}       \put(6,2.02){\line(1,0){3}}
                    \put(6,2.01){\line(1,0){3}}
                    \put(2.5,2){\line(1,0){3}}         \put(6,2){\line(1,0){3}}
                    \put(2.5,5){\line(1,0){3}}         \put(6,5){\line(1,0){3}}
           \end{picture}
     \caption{} \label{pic4}
     \end{figure}
          \\
       We want to show  that equation (\ref{equation sieben})  holds for all basis elements $T \in  {\cal S}_n(X)$.
       We have  that
        \begin{align*}
              ( \partial_{n+1} \circ \Theta_n )(T) &  =
              \partial_{n+1} \left( \sum_{ \vec{z} \in \Upsilon^{n} } \ \ v_{\vec{z}} \cdot  G_{\vec{z}}(T) \right) \\
           &  =  \sum_{j=1}^{n+1} (-1)^{j+1}   \sum_{ \vec{z} \in \Upsilon^{n}}   v_{\vec{z}} \cdot
                 \left( \: a \cdot  \left\langle \: G_{\vec{z}}(T) \:  \right\rangle_{0,j} + b \: \cdot                                 \left\langle \: G_{\vec{z}} (T) \:  \right\rangle_{1,j} \right) .
        \end{align*}
      In the beginning  let  us consider the special case   $j := n+1$.
      For   $\vec{z} := ( 0,0,0, \ldots , 0 ) \in \Upsilon^{n}$   we compute
             $\left\langle G_{\vec{z}} (T) \right\rangle_{1,n+1}$   for
            ($ x_1 , x_2, \ldots , x_n )  \in  { \bf I }^{n} $:
      \begin{align*}
             \left\langle G_{\vec{z}} (T) \right\rangle_{1,n+1} ( x_1 , x_2 , \ldots , x_n )
             &  =  G_{\vec{z}}(T)( x_1 , x_2,  \ldots , x_n ,1 )   \\
                =   T ( \eta_{0}(x_1,1), \eta_{0}(x_2,1),    \ldots ,   \eta_{0}(x_n,1) )
             &  =  T ( x_1 , x_2 ,   \ldots ,   x_n  ) .
      \end{align*}
       Hence       $b \cdot \left\langle G_{\vec{z}} (T) \right\rangle_{1,n+1} = b \cdot T$.
      (We are just computing the top of the squares, see again the previous figure!)
      We shall see that for the other    $\vec{z}  \in \Upsilon^{n}$   the corresponding elements of
             $\left\langle  G_{\vec{z}} (T) \right\rangle_{1,n+1}$  cancel pairwise. 
                
       For a fixed   $k  \in \{ 1,2, \ldots ,  n \}$   let \\
             $\vec{\lambda} := ( z_1, z_2, \ldots , z_{k-1} , 1 , z_{k+1} , \ldots , z_n ),
             \vec{\zeta} := ( z_1, z_2, \ldots , z_{k-1} , 2 , z_{k+1} , \ldots , z_n ) \in \Upsilon^{n}  $.
             For an arbitrary element ($ x_1 , \ldots , x_k , \ldots , x_n )   \in    { \bf I }^{n}$  we get
    \begin{align*}
           &  \left\langle G_{\vec{\lambda}} (T) \right\rangle_{1,n+1} ( x_1 , \ldots , x_k , \ldots , x_n )
                 =   G_{\vec{\lambda}} (T) ( x_1 ,  \ldots , x_k , \ldots , x_n , 1 )  \\
              &  =   T ( \ldots  \ldots , \eta_{1}(x_k,1) , \ldots    \ldots  )  \
                 =   T ( \ldots  \ldots , 1 , \ldots  \ldots  )  \\
              &  =   T (    \ldots   \ldots , \eta_{2}(x_k,1) , \ldots  \ldots  )  \
                 =   G_{\vec{\zeta}} (T) ( x_1 ,  \ldots ,  x_k , \ldots ,x_n , 1 )  \\
              &  =   \left\langle  G_{\vec{\zeta}} (T) \right\rangle_{1,n+1} ( x_1, \ldots , x_k, \ldots , x_n ) .
    \end{align*}
     Because  \ $v_{\vec{\lambda}}   \cdot  v_{\vec{\zeta}}  = -1$  it follows that
     $$ b \cdot v_{\vec{\lambda}} \cdot \left\langle G_{\vec{\lambda}} (T) \right\rangle_{1,n+1}  +  b
        \cdot v_{\vec{\zeta}} \cdot \left\langle  G_{\vec{\zeta}} (T)  \right\rangle_{1,n+1} = 0 .   $$
    (Now we compute the bottom of the squares, see again the previous figure.) We still have  $j = n+1$.
              We get for all  ($ x_1 , \ldots , x_n ) \in  { \bf I }^{n}$
              and  all   $\vec{z} =  ( z_1 , \ldots , z_n ) \in   \Upsilon^{n}$:
    \begin{align*}
          \left\langle G_{\vec{z}}(T) \right\rangle_{0,n+1} ( x_1 , x_2 , \ldots , x_n  ) &  =
          G_{\vec{z}}(T) \: ( x_1 , x_2 , \ldots , x_n , 0 )   \\
           =    T ( \eta_{z_1}(x_1,0), \eta_{z_2}(x_2,0),    \ldots ,   \eta_{z_n}(x_n,0) )
          &  =  T (  t_1 , t_2 ,   \ldots , t_n ) ,
    \end{align*}
        with   \quad \quad
        $ t_i :=  \begin{cases}    \frac{1}{3} \cdot  x_i   & \text{  if  }            z_i  = 0 ,  \\
                                   \frac{1}{3} \cdot (2 - x_i )  & \text{  if  }       z_i  = 1 ,  \\
                                   \frac{1}{3} \cdot  (2 + x_i ) & \text{  if  }       z_i  = 2  ,
                                              \ \text{ for   all} \  i = 1, 2, \ldots , n \, .  \\
                     \end{cases}  $  \\   \\
    We define  $\vec{e} := ( e_1 , e_2 ,   \ldots , e_n  )$, $\vec{v} := ( v_1 , v_2 , \ldots , v_n ) $,
           by setting for all $ i \in  \{1,2, \ldots n \}$:
           $$   e_i  :=  \begin{cases}   0  & \text{  if  }          z_i  = 0   \\
                                       2  & \text{  if  }             z_i \in \{ 1,2 \} \ ,    \\
                       \end{cases}   \qquad
              v_i  :=  \begin{cases}   1   & \text{  if  }        z_i \in  \{ 0,2 \}\
                                     ( { \rm hence} \ e_i \in \{ 0, 2 \} )    \\
                                       -1  & \text{  if  }        z_i = 1   ,
                                     ( { \rm hence} \ e_i = 2 ) \ .      \\
                       \end{cases}   $$
     We have  $ \vec{e} \in  {\cal E}^{n}$   and   $ \vec{v} \in  {\cal V}_{\vec{e},n}$.
     With  a few calculations it is easy to see that
       $$  v_{\vec{z}} \cdot \left\langle  G_{\vec{z}}(T)  \right\rangle_{0,n+1} \ = \
                 ( \prod_{i=1}^{n}  v_i ) \cdot \left\| T \right\|_{\frac{1}{3} , \vec{e} , \vec{v} } .  $$
      We compare this with the definition of  ${\cal SD}_n (T)$ in  (\ref{definition SD}).    We get
      $$     \sum_{ \vec{z} \in \Upsilon^{n}} \ a \cdot v_{\vec{z}} \cdot
                \left\langle  G_{\vec{z}}(T) \right\rangle_{0,n+1} \ = \ - a \cdot {\cal SD}_n (T).  $$
      All in all for the fixed  $j = n+1$ follows
       \begin{align}
                 \sum_{ \vec{z} \in \Upsilon^{n}}   v_{\vec{z}} \cdot
                 \left[ a \cdot  \left\langle G_{\vec{z}}(T)  \right\rangle_{0,n+1} + b  \cdot                                           \left\langle G_{\vec{z}} (T)  \right\rangle_{1,n+1}  \right]
                  =  b \cdot T  -  a  \cdot  {\cal SD}_n (T) .
       \end{align}
    Now let  $j$ be an element of  $ \{ 1, 2, 3, \ldots , n \}$. Let \\
             $\vec{\zeta_0} := ( z_1, z_2, \ldots , z_{j-1} , 0 , z_{j+1} , \ldots , z_n ), \
             \vec{\zeta_1} := ( z_1, z_2, \ldots , z_{j-1} , 1 , z_{j+1} , \ldots , z_n ) \in \Upsilon^{n}$,
             and  for all points  $( x_1 , x_2 , \, \ldots, \, x_n ) \in { \bf I }^{n}$  we get
    \begin{align*}
            &  \left\langle  G_{\vec{\zeta_0}} (T) \right\rangle_{1,j} ( x_1 , x_2 , \ \ldots, \, x_n )
              =   G_{\vec{\zeta_0}} (T) ( x_1 , x_2, \, \ldots, \, x_{j-1} , 1 , x_j , x_{j+1}, \, \ldots, \, x_n ) \\
            & =  T \left( \eta_{z_1}(x_1,x_{n}),     \ldots , \eta_{z_{j-1}}(x_{j-1},x_{n} ),  \
            \eta_{0}(1,x_n), \eta_{z_{j+1}}(x_{j},x_{n}), \ldots , \eta_{z_n}(x_{n-1},x_{n}) \, \right) \\
            & =   T \left( \eta_{z_1}(x_1,x_{n}), \ \ldots, \ \frac{1}{3 - 2 \cdot x_n}, \ \ldots \ \right)
              =  T \left( \eta_{z_1}(x_1,x_{n}), \, \ldots, \, \eta_{1}(1,x_n), \, \ldots \, \right)  \\
            & =   G_{\vec{\zeta_1}} (T) ( x_1 , x_2 , \ \ldots, \ x_{j-1} , 1 , x_j , x_{j+1} , \ \ldots , \ x_n )
              =  \left\langle G_{\vec{\zeta_1}} (T) \right\rangle_{1,j} ( x_1 , x_2 , \ \ldots , \ x_n )  .
     \end{align*}
  Thus,
$$b \cdot v_{\vec{\zeta_0}} \cdot \left\langle G_{\vec{\zeta_0}} (T) \right\rangle_{1,j}
        + \ b \cdot v_{\vec{\zeta_1}} \cdot \left\langle G_{\vec{\zeta_1}} (T) \right\rangle_{1,j}  =  0.$$(See the previous figure: the right hand side of  $T(\eta_0)$
cancels the right hand side of $-T(\eta_1)$.) With the same fixed
number   $j$   we calculate with the tuples
    $$\vec{\zeta_1} = ( z_1, z_2, \ldots , z_{j-1} , 1 , z_{j+1} , \ldots , z_n ),
           \vec{\zeta_2} := ( z_1, z_2, \ldots , z_{j-1} , 2 , z_{j+1} , \ldots , z_n ) \in
           \Upsilon^{n} $$
to obtain
\begin{align*}
&\left\langle G_{\vec{\zeta_1}} (T)  \right\rangle_{0,j} ( x_1 , x_2
, \ldots , x_n )\\ &=   G_{\vec{\zeta_1}} (T) ( x_1 , x_2 , \ldots ,
x_{j-1} , 0 , x_j, x_{j+1} , \ldots , x_n)\\
&= T (  \eta_{z_1}(x_1,x_{n}),     \ldots   \, \ldots ,
\eta_{1}(0,x_n)   ,\ldots    \ldots , \eta_{z_n}(x_{n-1},x_{n})  )\\
& = \ T ( \ldots  \ldots , t_j ,\ldots \ldots  ) \qquad
\text{ with }  t_j  := \begin{cases} \frac{2}{3 - 2 \cdot x_n} & \text{ if } x_n  \in  \left[ 0, \frac{1}{2} \right] , \\
                                    1  & \text{  if  }  x_n  \in  \left[ \frac{1}{2}, 1 \right]
                                    ,
                    \end{cases}\\
&= T (    \ldots   \ldots , \eta_{2}(0,x_n) ,\ldots \ldots )\\
&=  G_{\vec{\zeta_2}} (T) ( x_1 , x_2 , \ldots , x_{j-1} , 0 , x_j,
x_{j+1} , \ldots , x_n )\\
&= \left\langle G_{\vec{\zeta_2}} (T)  \right\rangle_{0,j}( x_1 ,
x_2 , \ldots , x_n ).
\end{align*}
Hence,
        $$a \cdot v_{\vec{\zeta_1}} \cdot \left\langle G_{\vec{\zeta_1}} (T) \right\rangle_{0,j}
        + a \cdot v_{\vec{\zeta_2}} \cdot \left\langle G_{\vec{\zeta_2}} (T) \right\rangle_{0,j}  =  0. $$
(See the previous figure again: The left side of $-T(\eta_1)$ cancels the left side of $T(\eta_2)$.)  \\
Now we take again  $\vec{\zeta_0} = ( z_1, z_2, \ldots , z_{j-1} , 0
, z_{j+1} , \ldots , z_n ) \in \Upsilon^{n}$,  define
$$\vec{\mu}=( z_1, z_2, \ldots , z_{j-1} , z_{j+1} , \ldots , z_n )
\in \Upsilon^{n-1},$$           and we get for ($ x_1 , \ldots ,
x_{j-1}, x_{j}, \ldots , x_n ) \in  { \bf I }^{n}$:
\begin{align*}
& \left\langle  G_{\vec{\zeta_0}} (T) \right\rangle_{0,j}  ( x_1,
\ldots , x_{j-1}, x_{j}, \ldots , x_n )\\
& = G_{\vec{\zeta_0}}(T) ( x_1 , \ldots , x_{j-1} , 0 , x_{j} ,
\ldots , x_n )\\
&= T ( \eta_{z_1}(x_1,x_{n}), \ldots ,
\eta_{z_{j-1}}(x_{j-1},x_{n}), \eta_{0}(0,x_{n})  ,
\eta_{z_{j+1}}(x_{j},x_{n}),  \ldots , \eta_{z_n}(x_{n-1},x_{n}))\\
&=  T (  \ldots  \ldots ,\eta_{z_{j-1}}(x_{j-1},x_{n}), 0 ,  \
               \eta_{z_{j+1}}(x_{j},x_{n}),  \ldots  \ldots)\\
&=\left\langle  T  \right\rangle_{0,j} ( \eta_{z_1}(x_1,x_{n}),
\ldots  \ldots , \ \eta_{z_{j-1}}(x_{j-1},x_{n}),
\eta_{z_{j+1}}(x_{j},x_{n}), \ldots \ldots) \\
&= G_{\vec{\mu}} \left( \left\langle  T   \right\rangle_{0,j}
\right)      ( x_1, \ldots , x_{j-1}, x_{j}, \ldots , x_n ).
\end{align*}
Hence $$a \cdot \left\langle  G_{\vec{\zeta_0}} (T)
\right\rangle_{0,j} =  a \cdot  G_{\vec{\mu}} \, \left( \left\langle
T \right\rangle_{0,j}         \right).$$ And, last but not least, we
can show in the same way  that
      $$b \cdot \left\langle  G_{\vec{\zeta_2}} (T) \right\rangle_{1,j}
                     \ = \  b \cdot  G_{\vec{\mu}} \left( \, \left\langle  T  \right\rangle_{1,j} \,
                     \right) $$
for $\vec{\zeta_2} = ( z_1, z_2, \ldots , z_{j-1} , 2 , z_{j+1} ,
\ldots , z_n ) \in \Upsilon^{n}$.

   Now we have collected all the needed  facts to confirm for every
         $n  \in \mathbbm{N}$   the equation   \\
   \centerline{  $( \partial_{n+1} \circ \Theta_n )(T) =  (-1)^{n+2}  \cdot
                 ( b \cdot Id - a \cdot {\cal SD}_n )(T) + ( \Theta_{n-1} \circ  \partial_{n}) (T)$.  }    \\
    Therefore, if we use a chain  {\sf u} instead of  $T$, the equation (\ref{equation sieben})
            is proved, because all used maps are linear.  It is  a trivial consequence  that for a cycle
            {\sf u} (i.e.  $ \partial_{n }({\sf u}) = 0$)  the equation
            $ [ a \cdot {\cal SD}_n $({\sf u})$   ]_{\sim} = [ b \cdot {\sf u} ]_{\sim}$
            follows   on the level of   homology classes,
            because  $b \cdot {\sf u} - a \cdot {\cal SD}_n ({\sf u})$ \, is in the image of  $\partial_{n+1}$.

The next step is to show  that for a  cycle {\sf u} the  equation $[
b \cdot \: {\cal SD}_n ({\sf u})   ]_{\sim} = [ a \cdot {\sf u}
]_{\sim}$ also holds  on the level of  homology classes. Looking  at
the  previous proof  this seems obvious, and we shall not explain it
in all details. The proof is nearly the same, we only have to modify
it by  `turning it upside down'. Instead of using the three
auxiliary functions    $\eta_0 ,   \eta_1 ,   \eta_2   , $ we  need
three others  $\widetilde{\eta_0},   \widetilde{\eta_1},
\widetilde{\eta_2} : { \bf I }^{2} \rightarrow { \bf I } $.

  For  $x,y\in [0,1]$ define:
\begin{align*}
\widetilde{\eta_0}(x,y) &:=\frac{x}{1+2 \cdot y},\\
 \widetilde{\eta_1}(x,y)   & := \begin{cases}  \frac{2 - x}{1 + 2 \cdot y}    & \text{  for }
                                                        y  \geq  \frac{1}{2} -  \frac{1}{2} \cdot x  \quad ,    \\
                                                      1 & \text{  else }  \quad  ,
                                      \end{cases}\\
\widetilde{\eta_2}(x,y)   &  := \begin{cases}  \frac{2 + x}{1 + 2
\cdot y}    & \text{  for }  y  \geq  \frac{1}{2} +  \frac{1}{2} \cdot x \quad ,     \\
                                                      1 & \text{  else } \quad .
                                      \end{cases}
\end{align*}
Then   $\widetilde{\eta_0} , \widetilde{\eta_1} , \widetilde{\eta_2}
$ are continuous. For a fixed tuple    $\vec{z} = ( z_1, z_2, \ldots
, z_n ) \in \Upsilon^{n}$ and for    $T \in  {\cal S}_{n}(X) $ let
us define the map $\widetilde{G_{\vec{z}}}(T)  : { \bf I }^{n+1}
\rightarrow X$ by setting for all   $( x_1, \ldots , x_n, x_{n+1} )
\in { \bf I }^{n+1}$:
     $$\widetilde{G_{\vec{z}}}(T)  ( x_1,  \ldots , x_n,  x_{n+1}  ) :=
                 T \left( \widetilde{\eta_{z_1}}(x_1,x_{n+1}), \widetilde{\eta_{z_2}}(x_2,x_{n+1}),    \ldots ,                               \widetilde{\eta_{z_n}}(x_n,x_{n+1}) \right)  .   $$
      Thus, $\widetilde{G_{\vec{z}}}(T) \in   {\cal S}_{n+1}(X) $. \ Let  for
                 $\vec{z} = ( z_1, z_2, \ldots , z_n ) \in \Upsilon^{n}$ the sign
                 $ v_{\vec{z}} := (-1)^{\sum_{i=1}^{n} z_i }$ as before, and finally  define
      \begin{align*}
                \widetilde{\Theta_n} (T)  :=  \sum_{ \vec{z} \in \Upsilon^{n} }  \quad   v_{\vec{z}} \cdot                              \widetilde{G_{\vec{z}}}(T) .
      \end{align*}
    By a similar calculation as for the proof of equation (\ref{equation sieben}) we show:
   \begin{align}
           ( \partial_{n+1} \circ \widetilde{\Theta_n}  )(T) =  (-1)^{n+2}  \cdot
           ( a \cdot Id - b \cdot {\cal SD}_n ) (T) + ( \widetilde{\Theta_{n-1}} \circ \  \partial_{n} ) (T) ,
   \end{align}
          and by using a cycle {\sf u} instead of the map $T$, it leads directly to the desired formula
    $$  [ b   \cdot   {\cal SD}_n ({\sf u}) ]_{\sim} = \ [a \cdot {\sf u} ]_{\sim} .  $$
     Now let   for all    $k  \in  \mathbbm{N},   n \! \in \! \mathbbm{N}_0$     and
              all chains   {\sf u}  $\in  {\cal K}_{ n} (X)$  \\
      \centerline{        ${\cal SD}_n^{(k)}$ ({\sf u}) :=
                         (${\cal SD}_n \circ {\cal SD}_n \circ   \ldots   \circ {\cal SD}_n$)
                         ({\sf u}),  \ (with   $k$   factors ${\cal SD}_n$).  }
    \begin{lemma}     \label{lemma vier}   \quad
          For all $k \in \mathbbm{N} , n \in \mathbbm{N}_0$,
                and  all \ {\sf u} $\in \text{kernel} (\partial_{n})$  we have  \\
    \centerline{     $ [  a^{k}\, \cdot {\cal SD}_n^{ (k)} ( {\sf u} ) ]_{\sim}
                     =   [ b^{k}  \cdot {\sf u} ]_{\sim}   $ \quad  and   \quad
                     $ [  b^{k}   \cdot  {\cal SD}_n^{ (k)} ( {\sf u} ) ]_{\sim}
                     =   [ a^{k} \cdot  {\sf u} ]_{\sim}  $ .    }
    \end{lemma}
   \begin{proof} \quad   We  prove the first equation by induction on $ k$.  Note that, if
           $   \partial_{n }$ ({\sf u}) = 0, also  $\partial_{n } ({\cal SD}_{n}^{\: (k)}(${\sf u})) = 0
           (because ${\cal SD}_{n}$ commutes with the boundary   operator), and note that ${\cal SD}_{n}$
           is a linear map.  Assume  for some  $k \in  \mathbbm{N}$ for an
           {\sf u} $\in \text{ kernel} (\partial_{n}):  \quad
           [ a^{k}  \cdot {\cal SD}_n^{ (k)}$ ({\sf u})$ ]_{\sim} =  [ b^{k} \cdot ${\sf u}$ ]_{\sim} $. \
           Further let  {\sf w} := $a^{k}  \cdot  {\cal SD}_n^{ (k)} $({\sf u}),
           hence  {\sf w} is a cycle, too.  Thus  we  get
    \begin{align*}
               [  b^{k+1} \cdot {\sf u} ]_{\sim}
             & =  [ b \cdot b^{k} \cdot  {\sf u} ]_{\sim}
               = [ b \cdot a^{k}   \cdot {\cal SD}_n^{ (k)} ({\sf u}) ]_{\sim}
               =  [   b   \cdot   {\sf w}   ]_{\sim}
               =   [   a   \cdot    {\cal SD}_n ({\sf w})   ]_{\sim}  \\
             & =  [ a \cdot a^{k}   \cdot {\cal SD}_n   ( {\cal SD}_n^{\: (k)} ({\sf u}) \, ) ]_{\sim}
               =   [ a^{k+1}   \cdot {\cal SD}_n^{\: (k+1)} ({\sf u}) ]_{\sim}  .
     \end{align*}
     This proves the first equation of the lemma for all $k  \in  \mathbbm{N}$.
    \end{proof}
    \begin{lemma}    \label{lemma fuenf}   \quad
             For all weights    $\vec{m} = \left( a ,b \right) \in  {\cal R}^{2}$,
                       $ (a,b) $   has the property   ${\cal {NCD}}$,   and for all
                       $k \! \in \! \mathbbm{N}$    there is an element   $r_k \in  {\cal R}$  such that the equation
                       $ r_k \cdot [ {\cal SD}_n^{\: (k)} ({\sf u}) ]_{\sim}  = [ {\sf u} ]_{\sim} $
                       holds for all  {\sf u} $\in \text{ kernel} (\partial_{n})$  and
                       $n \! \in \! \mathbbm{N}_0$   on the level of homology classes.
    \end{lemma}
    \begin{proof}  \quad
           The property ${\cal {NCD}}$  means that for all  $k  \in   \mathbbm{N}$   there are
                 $ x_k ,   y_k \in  {\cal R}$  such that  $ x_k \cdot a^{k} +  y_k  \cdot b^{k} = 1_{{\cal R}} $.
                 Now  set  $r_k  :=  x_k \cdot  b^{k}  +   y_k  \cdot    a^{k} $.
                 Take the previous   Lemma \ref{lemma vier}  and write
       \begin{align*}
                  [ {\sf u} ]_{\sim}  = [ 1_{{\cal R}} \cdot  {\sf u} ]_{\sim}
            &  =   [ (x_k \cdot  a^{k} + y_k  \cdot b^{k}) \cdot {\sf u} ]_{\sim}
               =    x_k \cdot  [ a^{k} \cdot {\sf u} ]_{\sim} + y_k  \cdot [ b^{k} \cdot {\sf u} ]_{\sim}  \\
            &  =    x_k \cdot  [  b^{k} \cdot {\cal SD}_n^{ (k)} ({\sf u}) ]_{\sim}  + \
                     y_k  \cdot [ a^{k} \cdot {\cal SD}_n^{ (k)} ({\sf u}) ]_{\sim}    \\
            &  =    ( x_k \cdot  b^{k} +  y_k  \cdot  a^{k} ) \cdot \
                      [ {\cal SD}_n^{ (k)} ({\sf u}) ]_{\sim}
               =    r_k \cdot  [ {\cal SD}_n^{ (k)} ({\sf u}) ]_{\sim} ,
        \end{align*}
               and the lemma is proved.
   \end{proof}
    Now let's return to the statement of Proposition \ref{Proposition Isomorphismus}.
    We still are proving  that for all $n \in  \mathbbm{N}$ the canonical inclusions
    ${\cal K}_{ n}( X,{\cal U} )  \stackrel{j}{\hookrightarrow}  {\cal K}_{ n}( X )$
          induce isomorphisms  \\  \centerline{  $j_\ast: \ _{\vec{m}}{\cal H}_{n}( X, \: {\cal U})
          \stackrel{\cong }{\longrightarrow} \, _{\vec{m}}{\cal H}_{n}( X )$. }

     Note that, for a  continuous  $T: { \bf I }^{n} \rightarrow X $, the image
          $T({ \bf I }^{n}$) is compact in $X$, and for the  given open covering
          $ \{ Int(U_i)  \: \mid \: i \in  \Im  \}$  of $X$  a finite subset is sufficient  to cover
          $T({ \bf I }^{n}$).   Further note that the diameters of the $3^{n}$ elements
          of the  chain   $ {\cal SD}_n (T)({ \bf I }^{n})$   decrease to a third,  compared with the
          diameter of   $T({ \bf I }^{n})$. Hence, by iterating  ${\cal SD}_n$,   there is a number
          $k_T \in \mathbbm{N}$    such that    ${\cal SD}_n^{(k_T)}(T) \in  {\cal K}_{ n}(X,{\cal U})$.
          And therefore for a chain {\sf u} $\in {\cal K}_{ n}(X)$ (which is a finite linear combination of
          some $ T's $), a number $k_{\sf u} \in \mathbbm{N}$  exists
          such that   ${\cal SD}_n^{(k_{\sf u})}$({\sf u}) $\in  {\cal K}_{ n}(X,{\cal U})$. \\

     Now let us look on the inclusion
         $j :   {\cal K}_{ n}( X,{\cal U} )  \stackrel{}{\hookrightarrow}  {\cal K}_{ n}( X )$,
         and the induced ${\cal R}$-module morphism
         $j_\ast : \, _{\vec{m}}{\cal H}_{n}( X, \: {\cal U})
         \stackrel{ }{\longrightarrow} \, _{\vec{m}}{\cal H}_{n}( X \: )$.
         We have to show that $j_\ast$ is an epimorphism and a monomorphism.  (compare  \cite[p.36]{Massey}).

     $j_\ast$   is an epimorphism: \\
        Let  {\tt z} $ \in \,  _{\vec{m}}{\cal H}_{n}( X \: )$.  Then it follows that
           there is a chain   {\sf u} $\in  \text{ kernel} (\partial_{n}) \subset {\cal K}_{ n}( X )$,
           with   $[ {\sf u} ]_{\sim} $  = {\tt z}.    Hence we can deduce that there is  a
           $k \in  \mathbbm{N}$  with
           ${\cal SD}_n^{(k)}$({\sf u}) $\in  {\cal K}_{ n}(X,{\cal U})$.  Take the factor
           $r_{k} \in   {\cal R}$  (see Lemma \ref{lemma fuenf})   and  write   \\
      \centerline { $ [ r_{k}   \cdot   {\cal SD}_n^{\: (k)} ({\sf u}) ]_{\sim}
           = [ {\sf u} ]_{\sim}  \in \, _{\vec{m}}{\cal H}_{n}(X)$ \ with
           $r_{k} \cdot  {\cal SD}_n^{(k)}({\sf u}) \in  {\cal K}_{ n}(X,{\cal U}) $. }
     $$  \text{Hence} \ \ j_\ast ( [ r_{k}   \cdot   {\cal SD}_n^{\: (k)} ({\sf u}) ]_{\sim} ) =
           [ j ( r_{k} \cdot {\cal SD}_n^{\: (k)} ({\sf u}) ) \: ]_{\sim}  =
           [ r_{k} \cdot {\cal SD}_n^{ (k)} ({\sf u}) ]_{\sim} = [ {\sf u} ]_{\sim}  = {\tt z}.  $$

      $j_\ast$   is a monomorphism: \\
      Let   {\tt x} $\in \,  _{\vec{m}}{\cal H}_{n}( X, {\cal U})$  with  $j_\ast$({\tt x}) = 0.
                      We must show that   {\tt x} = 0.

      For every   {\tt x}  $\in \, _{\vec{m}}{\cal H}_{n}( X, \: {\cal U}) $ exists a cycle
      {\sf v} $  \in {\cal K}_{ n}(X,{\cal U})$    with    $ [   {\sf v}   ]_{\sim}  $ = {\tt x}.                                       We must show that   {\sf v}   is a boundary, i.e. we have to show that
                there is a     ${\sf w} \in {\cal K}_{  {n+1}}(X,{\cal U})$
                with     $   \partial_{n+1} ({\sf w}) = {\sf v}$.    We have  \\
    \centerline{     $ [ j({\sf v}) ]_{\sim}   = j_\ast(  [ {\sf v} ]_{\sim})
                     = j_\ast( ${\tt x}$) = 0 \in \, _{\vec{m}}{\cal H}_{n}(X)$.        }    \\
     The assumption $j_\ast$({\tt x}) = $ 0 \in \, _{\vec{m}}{\cal H}_{n}(X)$     means that
                $j_\ast$({\tt x})   is the equivalence class   of a cycle  which is a boundary,
              i. e.   that there is a chain    $\widehat{{\sf w}} \in   {\cal K}_{ {n+1}}(X)$
              with   $   \partial_{n+1} (\widehat{{\sf w}}) = {\sf v} $.

     Choose a sufficient large number $k \in \mathbbm{N}$  such that ${\cal SD}_{n+1}^{\: (k)} (\widehat{{\sf w}})
             \in  {\cal K}_{ {n+1}}(X,{\cal U})$, and take the element   $r_k \in   {\cal R}$
              from  Lemma \ref{lemma fuenf}   and we have
               $[ r_k   \cdot   {\cal SD}_n^{\: (k)} ({\sf v}) ]_{\sim} = [ {\sf v} ]_{\sim}$ ,  (and note that
              ${\cal SD}_n^{\: (k)} ({\sf v}) \in {\cal K}_{ {n}}(X,{\cal U})$  by triviality).  Hence
     \begin{align*}
               \partial_{n+1} ( r_k   \cdot   {\cal SD}_{n+1}^{\: (k)} ( \widehat{{\sf w}}) )
                    =   r_k   \cdot   {\cal SD}_n^{\: (k)}  (   \partial_{n+1} ( \widehat{{\sf w}}) )
                    =   r_k   \cdot   {\cal SD}_n^{\: (k)}  ( {\sf v} ).
      \end{align*}
     So we conclude that   $\ r_k   \cdot   {\cal SD}_n^{\: (k)} ({\sf v}) $  is a boundary,
              therefore  $ [  r_k   \cdot   {\cal SD}_n^{\: (k)} ({\sf v}) ]_{\sim} = 0 $,  and  since
              $[ r_k   \cdot   {\cal SD}_n^{\: (k)} ({\sf v}) ]_{\sim} = [ {\sf v} ]_{\sim}  \,
              \in \, _{\vec{m}}{\cal H}_{n}( X, \: {\cal U})$  it follows  $[ {\sf v} ]_{\sim} = 0 = $ {\tt x}!
  \end{proof}
    \quad   Hence $ j_\ast : \, _{\vec{m}}{\cal H}_{n}( X, \: {\cal U})
       \stackrel{\cong }{\longrightarrow} \, _{\vec{m}}{\cal H}_{n}( X \: )$
       is an isomorphism, and the proof of   Proposition \ref{Proposition Isomorphismus}   is finished.
       This proposition leads directly to the excision axiom, see again  \cite[p.30,31]{Massey}.

\section{Computing Homology Groups.   Dividing by the Degenerate Maps}
In the last section we proved the excision axiom for a  `weight' \
$\vec{m} =  \left( m_{0}, m_{1}  \right)$   which has the property
${\cal {NCD}}$. In this way we constructed an extraordinary homology
theory. When you read the construction of this homology theory for
the first time it seems to be very difficult to  compute  the
homology modules of any space, except for a point. But fortunately
there is an old (1968) paper \cite{BuCoFl}, which helps us by using
the ordinary singular homology theory.
\begin{theorem}   \label{Theorem vier}
For  abelian groups ${\cal A}$ and all  $n \in  \mathbbm{N}_0$  and
all pairs of finite {\it CW}-complexes  $(X,B)$  let  $_{\cal
S}H_{n} [(X,B);{\cal A}]$ be the  $n^{th}$  ordinary singular
homology group  $($with coefficient  group  $_{\cal S}H_{0}
[(point)]  \cong {\cal A})$. Let  ${\cal R} :=  \mathbbm{Z}$. Let  $
\vec{m} =  \left( a , b \right) \in \mathbbm{Z}^{2}$ with
gcd$\{a,b\}=1$. Then we have  for each pair of finite {\it
CW}-complexes $(X,B)$ and for all $n \in \mathbbm{N}_0$ {\rm :}\\
If   $\{a,b\}  =  \{1,-1\}$ then there is an isomorphism $$
_{\vec{m}}{\cal H}_{n}(X,B) \cong
                           \sum_{k = 0}^{n} \ _{\cal S}H_{k} [(X,B);\mathbbm{Z}]  . $$
If   $ \{a,b\}  \neq   \{1,-1\}$, so that  the index   $\sigma = a+b
\neq 0$,  then
            \[     _{\vec{m}}{\cal H}_{n}(X,B)     \cong
                 \begin{cases}
                        \sum_{k \in \{ 0,1,2, \ldots \} \, \wedge \, 2k \leq n} \ _{ {\cal S} }H_{2k}
                        [(X,B);\mathbbm{Z}_{\sigma}]   \  &  \text{   if   } n \text{   is even },  \\
                        \sum_{k \in \{ 0,1,2, \ldots \} \, \wedge \, 2k+1 \leq n} \  _{ {\cal S} }H_{2k+1}
                        [(X,B);\mathbbm{Z}_{\sigma}]  \  &  \text{   if   } n \text{   is odd}.
                 \end{cases}
            \]
   \end{theorem}
   \begin{proof} \ See \cite{BuCoFl}, and use the homology groups of a point, computed in section 3.
   \end{proof}
Recall that in section \ref{section boundary operator} we calculated
the homology groups of a point for an arbitrary ${\cal R}$, and note
that `our' homology theory   $ _{\vec{m}}{\cal H}_{n}$   differs
from the usual singular homology theory. But, as we announced in the
abstract,  we can divide the chain modules  ${\cal K}_{ n}(X)$  by
suitable submodules, and  in the case of  $\vec{m} = \left( m_{0},
m_{1} \right)$  where $(m_{0}, m_{1})$ has the property ${\cal
{NCD}}$, we shall obtain  the usual singular homology  theory  with
the coefficient module  ${\cal R} / ({\sigma}{\cal R})$. Compare
\cite[p.12,13]{Massey}, or \cite[p.236 ff]{Mac Lane}, where this
process is called a {\it normalization}.

As always, a new part begins with  definitions, see the definitions
\ref{definition sechs} - \ref{definition zehn}.
\begin{lemma}   \label{lemma sechs} \quad
For all weights   $\vec{m} =  \left( m_{0}, m_{1} , \ldots , m_{L}
\right) \in {\cal R}^{L+1}$ \ $($hence its index is $ \sigma =
\sum_{i=0}^{L} m_i ) , $ for $ n \in  \mathbbm{N}$ the boundary
operator  $ \partial_{n} :   {\cal K}_{ n}(X) \; \longrightarrow {\cal K}_{ n-1}(X)$ yields a map \\
\centerline  { $ \partial_{n   \mid \Gamma_{\sigma,   n}(X)} :   \
 \Gamma_{\sigma,   n}(X) \longrightarrow  \Gamma_{\sigma,   n-1}(X)$.  }
   \end{lemma}
   \begin{proof}   \quad
           Let  {\sf u}  $\in \Gamma_{\sigma, n}(X)$.  We know that   \\
           \centerline {  {\sf u} = {\sf u}$_1$ + {\sf u}$_2$, \  with   {\sf u}$_1  \in$
                         $\Ideal _{\sigma, n}(X)$, {\sf u}$_2  \in {\cal K}_{{\cal D}, n}(X) $.    }  \\
           Hence it follows that $ \partial_{n}$({\sf u}$_1) \in \Ideal _{\sigma,  n-1}(X)$,
                 since $ \partial_{n}$ is linear.  And we have  {\sf u}$_2 = \sum_{k=1}^{p} r_k \cdot T_k $,
                 and all the   $T_k's$  are  degenerate. Take $T$ := $T_k$, and assume that $T$ is degenerate at the
                 $\widehat{j}^{th}$  component,  $\widehat{j}   \in \{ 1,2,  \ldots , n \}$, i.e.
                 for all  $y,z  \in [ 0,1 ]$  we have \ (see Definition \ref{definition sechs}) \\
      \centerline{   \ $T ( x_1, x_2, \ldots , x_{\: \widehat{j}-1} , y , x_{\: \widehat{j}+1} , \ldots , x_n )                      = \ T ( x_1, x_2, \ldots , x_{\: \widehat{j}-1} , z , x_{\: \widehat{j}+1} , \ldots , x_n )$. } \\ \\
        By the definition of $\partial_{n} (T)$ it follows that
             $$  \partial_{n}  (T)
             \ = \ \sum_{j \in \{1,2,\ldots n \} \wedge j \neq \widehat{j} } (-1)^{j+1} \cdot \sum_{i = 0}^{L} m_{i}
             \cdot \left\langle  T\right\rangle_{\: n,\;  i,\; j \;}  \quad + \quad (-1)^{ \widehat{j} +1}
             \cdot  \sum_{i = 0}^{L} m_{i} \cdot  \left\langle  T \right\rangle_{\: n,\;  i,\;  \widehat{j} \;} .   $$
    The first summand is a linear combination of degenerate maps.  We compute the second summand. For a point
             ($ x_1, x_2, \ldots , x_{n-1} $) $\in$ ${ \bf I }^{n-1}$   we have:
    \begin{align*}
       & \sum_{i = 0}^{L} m_{i} \cdot \left\langle T \right\rangle_{\: n,\;  i,\;  \widehat{j} \;}
              ( x_1, x_2, \ldots ,  x_{\: \widehat{j}-1} , x_{\: \widehat{j}} , \ldots ,  x_{n-1}  ) \\
       &  =   \sum_{i = 0}^{L} m_{i} \cdot \  T \left( x_1, x_2, \ldots , x_{\: \widehat{j}-1} ,
                          \frac{i}{L} , x_{\: \widehat{j}} , \ldots , x_{n-1} \right)  \\
       &  =   T ( x_1, \ldots , x_{\: \widehat{j}-1} , * , x_{\: \widehat{j}} , \ldots , x_{n-1} )  \cdot
              \sum_{i = 0}^{L} m_{i} \\
       &  =   T ( x_1, \ldots , x_{\: \widehat{j}-1} , * , x_{\: \widehat{j}} , \ldots , x_{n-1} ) \cdot \sigma  .
    \end{align*}
        Thus, the second summand is  an element of $\Ideal _{\sigma,  n-1}(X)$.
  \end{proof}
         In Definition   \ref{definition zehn} we defined the quotient ${\cal R}$-module
              $  {\cal K}_{ n}(X,A)_{\sim  \Gamma, \sigma }$.   By the previous Lemma  \ref{lemma sechs}
              the boundary operators   $ (\partial_{n})_{n \geq 0} $    yield  a chain complex
     \[   \cdots \cdots  \  \labto{\partial_{n+1}}       {\cal K}_{ n}(X,A)_{\sim  \Gamma, \sigma }                                                   \labto{ \ \partial_{n} \ }   {\cal K}_{ n-1}(X,A)_{\sim  \Gamma, \sigma }
                              \labto{ \partial_{n-1}}   \ \cdots \  \labto{ \ \partial_0 \ } \{ 0 \}   .
     \]
    This leads to homology  ${\cal R}$-modules as usual,   \
              $  _{\vec{m}}{\cal H}_{n /\Gamma}(X,A) := \frac {kernel (\partial_{n})}
              {image (\partial_{n+1})}, \, \text{for} \ n \in  \mathbbm{N}_0$.

    Example:  In section 3 we calculated the homology groups for a one-point space  $\{p\}$.
              We had for   $\sigma =  0$   that
               $_{\vec{m}}{\cal H}_{n}(p)   \cong  {\cal R} $  for  all  $n  \in \mathbbm{N}_0$,   and for arbitrary
               indexes  $\sigma $   we got:
          \[  _{\vec{m}}{\cal H}_{n}(p)
                     \cong \begin{cases} \{ x \in  {\cal R} \, | \, \sigma \cdot x = 0 \}
                                                                & \text{  if  } n  \text{ is odd }  \\
                     {\cal R}/(\sigma \cdot {\cal R})           & \text{  if  } n  \text{ is even}  .
                     \end{cases}
          \]
     For the space $\{p\}$   and for $n \in \mathbbm{N}$   the single map
           $T: { \bf I }^{n} \rightarrow \{p\}$ is  degenerate, but  $T: { \bf I }^{0} \rightarrow \{p\}$  is not.
           Hence   $ \Gamma_{\sigma,0 }(p)$ =$\Ideal _{\sigma,0}(p) \cong  \sigma \cdot {\cal R} $,
           thus the generating chain complex  \quad
           $ _{\vec{m}}{\cal K}_*(p)\ = \ \cdots \ \labto{\partial_4} {\cal K}_{3}(p)
                  \labto{\partial_3} {\cal K}_{2}(p) \labto{\partial_2} {\cal K}_{1}(p)
                  \labto{\partial_1} {\cal K}_{0}(p) \labto{\partial_0} \{ 0 \}  ,    $
   \[ \text{ i. e. } \quad
    _{\vec{m}}{\cal K}_*(p)\ \cong  \  \cdots \  \labto{\partial_4}  {\cal R} \labto{\partial_3}  {\cal R}
                   \labto{\partial_2}  {\cal R}  \labto{\partial_1}  {\cal R}   \labto{\partial_0} \{ 0 \} ,
   \]
      turns, by dividing for each   $n \in \mathbbm{N}_0$   by   $\Gamma_{\sigma, n }(p)$,  into
    \begin{align*}  &  \cdots  \  \labto{\partial_4} {\cal K}_{3}(p)_{\sim  \Gamma, \sigma }
                                    \labto{\partial_3}  {\cal K}_{2}(p)_{\sim  \Gamma, \sigma }
                                  \labto{\partial_2}  {\cal K}_{1}(p)_{\sim  \Gamma, \sigma }
                                  \labto{\partial_1}  {\cal K}_{0}(p)_{\sim  \Gamma, \sigma }
                                  \labto{\partial_0}  \{ 0 \}  \\
            \cong \ &  \cdots  \  \labto{\partial_4}  \{ 0 \}  \labto{\partial_3}  \{ 0 \}
                                  \labto{\partial_2}  \{ 0 \}  \labto{\partial_1} {\cal R}/(\sigma \cdot {\cal R})
                                  \labto{\partial_0}  \{ 0 \}   .
    \end{align*}
      \[    \text{Hence it follows that}  \quad  \quad  _{\vec{m}}{\cal H}_{n /\Gamma}(p)
                    \cong   \begin{cases}
                             \{ 0 \}   \quad                             &   \text{ for }   n   \in \mathbbm{N}  \\
                                 {\cal R}/(\sigma \cdot {\cal R})  \quad &   \text{ for }   n   =   0 .
                            \end{cases}
      \]
    \begin{corollary}  \label{Korollar eins} \
            If we take a weight   $\vec{m} =  \left( a, b \right) \in {\cal R}^{2}$, and if   $ ( a , b ) $
                  has the property   ${\cal {NCD}}$, the homology theory
                  $ _{\vec{m}}{\cal H}_{/\Gamma} := (_{\vec{m}}{\cal H}_{n /\Gamma})_{ n \geq 0 }$
                  is isomorphic to the ordinary singular  homology  theory on all pairs of  finite
                  {\it CW}-complexes, and we have a coefficient module
                  $ _{\vec{m}}{\cal H}_{0 /\Gamma}(p) \cong {\cal R}/ (\sigma\cdot{\cal R}) $, with $\sigma = a + b $.
    \end{corollary}
    \begin{proof}  \quad
           As  we proved  above, the homology theory   $_{\vec{m}}{\cal H}$ fulfils all of the
                 Eilenberg-Steenrod axioms except one. That means the axioms of exactness,
                 homotopy and excision are satisfied.
                 By dividing   the chain modules  ${\cal K}_{ n}(X,A)$  by  $\Gamma_{\sigma, n }(X,A)$  and using the
                 boundary operator  $ \partial_{n} $  we get the  homology modules
                 $_{\vec{m}}{\cal H}_{n /\Gamma}(X,A)$, and the dimension axiom  will be added, while the other
                 three axioms remain. Thus, with the uniqueness theorem proved by  Eilenberg and Steenrod,
                 we have the  uniqueness of the homology groups for all finite
                 {\it CW}-complexes $(X,A)$. See  \cite[p.51 ff]{Hu}, or  \cite[p.100 ff]{Eil-Stee}.
    \end{proof}
    \begin{corollary} \
           The usual singular homology theory is a special case of the class which is developed here.
                 If we take  the weight $\vec{m} := \left( 1 , -1 \right)  \in  \mathbbm{Z}^{2}$,
                 the   homology theory   $_{\vec{m}}{\cal H}_{/\Gamma} $
                 is isomorphic to the usual  singular  homology theory on all pairs of  finite {\it CW}-complexes,
                 and the coefficient group is   $\mathbbm{Z} $.
    \end{corollary}
    \begin{proof} \quad
           By the previous Corollary  \ref{Korollar eins}. Or see for the last time \cite[p.11-37]{Massey}.
    \end{proof}

   \section{Final Suggestion}
         We cannot decide whether this  new homology theory has any  important application.
            Perhaps it might be an interesting tool for other mathematicians.
            It is easy to see one difficulty:
            The computation of the homology modules of the one-point space is very simple. But to do the same
            for other topological spaces might be more complicate (except for finite {\it CW}-complexes, see
            Theorem  \ref{Theorem vier}),  although the homotopy axiom and
            the excision axiom and the exactness axiom hold. There are not enough  $\{ 0 \}$'s  in the
            homology modules of  a point, e.g. for ${\cal R} :=  \mathbbm{Z}$ and $ \sigma \neq 0 $
            only every second is the trivial group $ \{ 0 \} $.
            So it would be an improvement for a better  application if we can increase the number of  $ \{ 0 \}$'s.

       If we consider for all   $ n \in  \mathbbm{N}_0 $   the canonical quotient map  $t_n $,
    $$t_n : {\cal K}_{n}(X) \longrightarrow  \frac{{\cal K}_{n}(X)}{\Gamma_{\sigma,n }(X)}
             \, ,  \quad  \text{then the following diagram commutes}: $$
         \[   \xymatrix{   \cdots \quad   \ar[r]^<>(0.5){   \partial_{3}}
          &  {\cal K}_{ 2}(X)   \ar[d]<0.0cm>^<>(0.4){t_2}
                                      \ar[r]^<>(0.5){   \partial_{2}}
          &  {\cal K}_{ 1}(X) \ar[d]<0.0cm>^<>(0.4){t_1}
                                      \ar[r]^<(0.4)>(0.5){   \partial_{1}}
          &  {\cal K}_{ 0}(X) \ar[d]<0.0cm>^<>(0.4){t_0}
                                      \ar[r]^<>(0.5){   \partial_{0}}
          &   \{ 0 \}                \ar[d]<0.0cm>^<>(0.4){0}    \\
         \cdots   \quad    \ar[r]_<>(0.5){   \partial_{3}}
          &  {\cal K}_{2}(X)_{\sim \Gamma, \sigma }   \ar[r]_<>(0.5){   \partial_{2}}
          &  {\cal K}_{1}(X)_{\sim \Gamma, \sigma }   \ar[r]_<>(0.5){   \partial_{1}}
          &  {\cal K}_{0}(X)_{\sim \Gamma, \sigma }   \ar[r]_<>(0.5){   \partial_{0}}
          &          \{ 0 \}                                                                  \\     }
        \]     \\
     Let   $\beta$   be a fixed element from the set   $\mathbbm{N}_0 \cup \{ \infty$ \}.
               Then we generate a chain complex
     \[   \cdots \cdots \    \labto{{\overline{\partial_{n+1}}}}  \overline{_{\beta }{\cal K}_{n}(X)}
                               \labto{{\overline{\partial_{n}}}}    \overline{_{\beta }{\cal K}_{n-1}(X)}
                               \labto{{\overline{\partial_{n-1}}}}   \ \cdots \
                             \labto{ \ {\overline{\partial_{1}}} \ }  \overline{_{\beta }{\cal K}_{0}(X)}
                             \labto{ \ {\overline{   \partial_{0}}}  \ } \{ 0 \}
     \]
            if we define  for all  $n \in   \mathbbm{N}_0$:
     \[        \overline{_{\beta }{\cal K}_{n}(X)}
                    :=   \begin{cases}
                              {\cal K}_{n}(X) & \quad  { \rm if } \quad  n  \geq \beta      \\
                              {\cal K}_{n}(X)_{\sim \Gamma, \sigma } & \quad { \rm if } \quad 0 \leq n < \beta .
                          \end{cases}
     \]
      If we have chosen a number $\beta \in \mathbbm{N}$, then let
      $\overline{\partial_{\beta }}: {\cal K}_{\beta }(X) \rightarrow {\cal K}_{\beta - 1 }(X)_{\sim \Gamma, \sigma }$
      be the $\beta^{th}$  boundary operator by defining
      $\overline{\partial_{\beta }} := \ \partial_{\beta } \circ t_\beta =
             t_{\beta -1} \circ \ \partial_{\beta }$.
      This means that for the special  cases  $\beta = 0$ and $\beta = \infty$ we have for all $n \in \mathbbm{N}_0$:
      $$   \overline{_{0}{\cal K}_{n}(X)}  =  {\cal K}_{n}(X)  \quad \text{ and } \quad
           \overline{_{\infty}{\cal K}_{n}(X)} = {\cal K}_{n}(X)_{\sim \Gamma, \sigma }, \quad \text{respectively}. $$        For arbitrary weights $ \vec{m} $ we define for  $n  \in   \mathbbm{N}_0$:
            $_{\vec{m},\beta}{\cal H}_{n}(X) := \frac{kernel(\overline{\partial_n})}
            {image(\overline{\partial_{n+1}})}$,
            and we get sequences
             $ _{\vec{m},\beta}{\cal H} := (_{\vec{m},\beta}{\cal H}_n)_{ n \geq 0 }$
             of homology modules; with the two we have developed  here as  special cases,
            i.e.   for  $\beta = 0$ and $\beta = \infty$ we get
            $$_{\vec{m}, 0 }{\cal H} = \, _{\vec{m}}{\cal H}  \quad   \text{and}  \quad
            _{\vec{m}, \infty }{\cal H}  = \, _{\vec{m}}{\cal H}_{/\Gamma}  .   $$                                              Example:   For the one-point space $\{ p \}$ (which is our favourite topological space obviously)
            and  ${\cal R} := \mathbbm{Z}$  and for the weight $\vec{m} := (1,4)$   and   $\beta := 7 $                            or   $\beta := 8$     we obtain for   $n \in   \mathbbm{N}_0$:
     \[  _{(1,4),7}{\cal H}_{n}(p)
             \cong
           \begin{cases}
             \mathbbm{Z}_5  \quad &     \text{  for }   \quad n \in   \{   0, 8, 10, 12, 14,   \ldots   \} \\
             \mathbbm{Z}   \quad  &     \text{ for  }   \quad   n = 7      \\                                                         \{ 0 \}      \quad  &     \text{  for }   \quad  n \in \{ 1, \ldots , 6, 9, 11, 13, 15, \ldots \} \, ,              \end{cases}
     \]
     \[   _{(1,4) ,8}{\cal H}_{n}(p)
             \cong
           \begin{cases}
             \mathbbm{Z}_5  \quad &  \text{ \rm for }   \quad n \in   \{   0 , 8 , 10 , 12 , 14 ,   \ldots   \} \\
              \{ 0 \}  \quad      &  \text{ \rm for }   \quad  n \in \{ 1, \ldots , 7, 9, 11, 13, 15, \ldots \} \, .
           \end{cases}
     \]       \\   \\
       Acknowledgements:  \quad
         The author likes to thank
         Dr. Bj\"orn R\"uffer, Prof. Dr. Eberhard Oeljeklaus, Dr. Nils  Th\"urey, Jan  Osmers, Walter Meyer,
         Prof. Dr. Rick Jardine and Prof. Dr. Hans-Eberhard Porst for interest and suggestions, discussions and
         technical help. Also we thank Prof. Dr. Ronald Brown for many hints and much patience, and
         Dr. Guentcho  Skordev, who supported us with attentive listening and clever remarks, and not to forget
         the unknown referee of the {\it Journal of Homotopy and Related Structures}
         for careful reading of the paper and suggesting many improvements.


\begin{thebibliography}{99}
\bibitem{BHS}
R.~Brown, P.~J. Higgins, and R.~Sivera.
\newblock {\em Nonabelian algebraic topology: filtered spaces, crossed
  complexes, cubical homotopy groupoids}.
\newblock EMS Tracts in Mathematics Vol. 15 (to appear 2011).
\bibitem{BuCoFl}  R.O. Burdick, P.E. Conner, E.E. Floyd, `Chain Theories and Their Derived Homology',
         \textit{Proc. Amer. Math. Soc.}, Vol. 19, No. 5 (1968), 1115-1118.
\bibitem{Eil-Stee} S. Eilenberg, N. Steenrod, {\it Foundations of Algebraic Topology},  Princeton, 1952.
\bibitem{Hatcher} A. Hatcher, {\it Algebraic Topology},   Cambridge University Press, 2001.
\bibitem{hilton-wylie}
P.~J. Hilton and S.~Wylie.
\newblock {\em Homology Theory: {A}n Introduction to Algebraic Topology}.
\newblock Cambridge University Press, New York, 1960.
\bibitem{Hu} S.T. Hu, {\it Homology Theory},  Holden-Day,  1966.
\bibitem{Ka-Mi-Mro} T. Kaczynski, K. Mischaikow, M. Mrozek, {\it Computational Homology},  Springer, 2004.
\bibitem{Mac Lane} S. Mac Lane, {\it Homology}, Springer, 1963, Fourth Printing 1994.
\bibitem{Massey} W.S. Massey, {\it Singular  Homology  Theory},     Springer,  1980.
\bibitem{Rotman} J. Rotman, {\it An Introduction to Algebraic Topology},   Springer, 1988.
\bibitem{Schoen} R. Sch\"on, `Acyclic Models and Excision', \textit{Proc. Amer. Math. Soc.},
         Vol. 59, No. 1 (1976),  167-168.
\bibitem{Spanier} E.H. Spanier, {\it Algebraic Topology},  Springer,  1966.
\bibitem{Vick} J. W. Vick, {\it Homology Theory},  Springer, 1994.
\end{thebibliography}
\end{document}